\documentclass{article}
\usepackage{arxiv}
\usepackage[utf8]{inputenc}
\usepackage[T1]{fontenc}
\usepackage{amsmath,amssymb,amsfonts,amsthm}
\usepackage{graphicx}
\usepackage{booktabs}
\usepackage{siunitx}
\usepackage{placeins}
\usepackage{tikz}
\usetikzlibrary{arrows.meta,positioning}
\usepackage{pgfplots}
\pgfplotsset{compat=1.16}
\usepackage{algorithm}
\usepackage{algpseudocode}
\usepackage{microtype}
\usepackage[numbers,sort&compress]{natbib}
\usepackage{hyperref}
\usepackage[capitalise,nameinlink]{cleveref}

\newtheorem{definition}{Definition}
\newtheorem{theorem}{Theorem}
\newtheorem{corollary}{Corollary}
\theoremstyle{remark}
\newtheorem{remark}{Remark}

\renewcommand{\shorttitle}{DQ-DLCT: Fast Algorithms and Pose Denoising}
\date{}
\title{Dual Quaternion Discrete Linear Canonical Transforms:
FFT-Based Fast Algorithms and Joint Rotation--Translation Pose Denoising}
\author{
Yi Wu \qquad Xuan Chen \qquad Kit Ian Kou\thanks{Corresponding author: Kit Ian Kou (\texttt{kikou@um.edu.mo}).}\\[0.4em]
\normalfont University of Macau\\
\texttt{koreywuyi@gmail.com} \quad \texttt{cx277787820@gmail.com}\\
\texttt{kikou@um.edu.mo}
}
\hypersetup{
  colorlinks=true,
  linkcolor=black,
  citecolor=black,
  urlcolor=black,
  pdftitle={Dual Quaternion Discrete Linear Canonical Transforms: FFT-Based Fast Algorithms and Joint Rotation--Translation Pose Denoising},
  pdfauthor={Yi Wu, Xuan Chen, Kit Ian Kou},
  pdfkeywords={dual quaternions; linear canonical transforms; fast algorithms; convolution theorems; rigid-body pose denoising}
}
\begin{document}
\maketitle
\begin{abstract}
Rigid-body pose signals couple rotation and translation, yet most
existing denoising pipelines treat them separately and analyse them
only in the Fourier domain. We show that this separation is
unnecessary and often suboptimal. We fuse dual quaternions, which
encode rotation and translation in a single singularity-free algebra,
with linear canonical transforms, which supply flexible
quadratic-phase analysis, to develop one- and two-dimensional dual
quaternion discrete linear canonical transforms (DQ-DLCTs) with an
explicit kernel ordering. We establish inverse and energy-preservation
properties, structured spatial and frequency convolution relations,
and, crucially, fast algorithms built on complex FFTs. These
algorithms make the DQ-DLCT practical: a two-stage chirp-adaptive
scheme then performs joint rotation--translation denoising in a single
transform domain. Numerical experiments confirm the accuracy and
efficiency of the fast transform and demonstrate strong denoising
performance on measured pose sequences and constructed
two-dimensional pose fields. The DQ-DLCT thus turns a heavy
algebraic construction into a deployable tool for pose signal
processing.
\end{abstract}


\section{Introduction}

Rigid-body pose signals carry rotation and translation as a coupled
pair, yet most existing denoising pipelines split the two apart and
analyse them only in the Fourier domain. This separation is
unnecessary, and it is often suboptimal. Dual quaternions encode the
coupling compactly and have been used in pose estimation, calibration,
and robotic kinematics
\cite{walker1991dual,daniilidis1999handeye,wang2024dual,
farias2024review,chandrasekar2025kinematic}; their algebraic role in
motion construction and factorization has also been examined in the
computational mathematics literature
\cite{li2019factorization,yaqub2026rational}. Together with methods
for approximating rigid-body motion \cite{vrablikova2024interpolation},
these developments support the view that pose data should be analysed
as a structured signal rather than as unrelated rotation and position
coordinates.

The linear canonical transform (LCT) generalizes Fourier analysis
through quadratic-phase kernels parameterized by a unimodular matrix
\cite{moshinsky1971canonical,healy2016linear}. Its discrete
realization has motivated work on sampling and numerical computation
\cite{pei2000closed,healy2009sampling,hennelly2005fast,
koc2008digital,koc2021operator}, and its canonical parameters offer a
flexible representation of chirp-like and nonstationary components
\cite{xia2000chirp,li2026synchrosqueezed}. For pose signals, however,
applying a scalar LCT componentwise does not by itself preserve the
algebraic coupling between rotation and translation.

Quaternion Fourier transforms supply spectral tools for
multicomponent signals \cite{sangwine1996colour,
ell2007hypercomplex,hitzer2021introduction}, and efficient
implementations can be built from complex FFTs
\cite{pei2001quaternionfft}. Since quaternion multiplication is
noncommutative, both kernel placement and factor order affect the
resulting convolution relations
\cite{bujack2017convolution,sfikas2025matrix}. These issues persist
for canonical transforms with quaternion-valued kernels. None of
these formulations, however, yields a discrete canonical transform in
the dual-quaternion algebra---the setting required for pose signals.

Closest to the present work, Bahri and Ashino developed a continuous
quaternion LCT \cite{bahri2019qlct}, while Urynbassarova, Teali, and
Zhang introduced a discrete quaternion LCT and studied its fast
computation and convolution properties \cite{bahri2021discrete}.
Further quaternion LCT convolution structures were developed in
\cite{hu2023convolution,urynbassarova2023offset}. Recent work has
also considered windowed reduced-biquaternion canonical transforms
\cite{yang2024reduced}, a quaternionic hyperbolic offset LCT
\cite{elhaoui2026hyperbolic}, and a discrete biquaternion LCT
\cite{ahmad2026discrete}. These constructions target different signal
domains or algebras: none supplies the nilpotent dual unit required
for a dual-quaternion pose representation.

A dual quaternion Fourier transform has been investigated for dual
quaternion valued signals \cite{kenwright2023dual}. The present work
instead develops a discrete linear canonical transform directly in the
dual-quaternion algebra, where kernel multiplication must be ordered
explicitly and the primal and dual components must remain coupled
through inversion, convolution, and fast computation. The resulting
transform is more than an algebraic curiosity: a
chirp--Fourier--chirp factorization reduces it to a fixed number of
ordinary complex FFTs, so the two-dimensional transform runs in
$\mathcal O(N_1N_2(\log N_1+\log N_2))$ time rather than the
$\mathcal O(N_1^2N_2^2)$ of direct evaluation. This is what makes the
construction deployable, and it is what allows us to use the transform
as a denoising tool: spectral concentration guides canonical-rate
selection \cite{baraniuk2001measuring}, while spectral subtraction and
shrinkage motivate the filtering step
\cite{boll1979,lim1979enhancement,donoho1994ideal}. The method is
evaluated on measured Mini Wheelbot pose sequences and on
two-dimensional pose fields constructed from its trajectories
\cite{hose2026wheelbot}.

The main contributions of this paper are as follows.
\begin{enumerate}
    \item One- and two-dimensional dual quaternion discrete linear
    canonical transforms (DQ-DLCTs) are constructed with an explicit
    kernel ordering, and their spatial and frequency convolution
    theory is established.

    \item Fast forward and inverse algorithms are developed based on a
    chirp--Fourier--chirp factorization and ordinary complex FFTs, and
    their accuracy and runtime are verified against direct
    computation.

    \item DQ-DLCT-based pose denoising is evaluated on measured
    one-dimensional sequences and constructed two-dimensional pose
    fields, demonstrating joint rotation--translation filtering in a
    single transform domain.
\end{enumerate}

The remainder of this paper is organized as follows.
Section~\ref{sec:prelim} provides the necessary background.
Section~\ref{sec:DQ-DLCT} develops the DQ-DLCT together with its fast implementation. Section~\ref{sec:numerical_properties} establishes the numerical properties of the two-sided transform.
Section~\ref{sec:convolution} presents the convolution theory, and
Section~\ref{sec:experiments} reports the numerical and
pose-denoising experiments. The advantages, limitations, and broader
application areas of the DQ-DLCT are discussed in
Section~\ref{sec:conclusion}.

\section{Preliminaries}\label{sec:prelim}

This section reviews the algebraic foundations of dual quaternions
and the linear canonical transform (LCT), and establishes the
notation and kernel identities used throughout the paper.

\subsection{Dual Quaternions}\label{subsec:dual_quaternion}

A dual quaternion is written as
$\widehat q = q_r+\varepsilon q_d$, with $q_r,q_d\in\mathbb H$,
where $\varepsilon$ is the dual unit satisfying
$\varepsilon^2=0$ and $\varepsilon\neq 0$. The dual unit
$\varepsilon$ commutes with all quaternion-valued quantities. The
set of dual quaternions is denoted by
\begin{equation*}
\mathbb{DH}
=
\left\{
q_r+\varepsilon q_d:
q_r,q_d\in\mathbb H
\right\}.
\end{equation*}

For the quaternion units $\mathbf{i}$, $\mathbf{j}$, and
$\mathbf{k}$, we have
$\mathbf{i}^2=\mathbf{j}^2=\mathbf{k}^2=\mathbf{i}\mathbf{j}\mathbf{k}=-1$
and
$\mathbf{i}\mathbf{j}=\mathbf{k}$, $\mathbf{j}\mathbf{i}=-\mathbf{k}$.
Hence quaternion multiplication is noncommutative
\cite{hitzer2021introduction}. This noncommutativity is the source
of the left--right asymmetry that distinguishes the transforms
defined in Section~\ref{sec:DQ-DLCT}.

We define the two complex subalgebras
\begin{equation*}
\mathbb C_{\mathbf{i}}=\{x+\mathbf{i}y:x,y\in\mathbb R\},
\qquad
\mathbb C_{\mathbf{j}}
=\{x+\mathbf{j}y:x,y\in\mathbb R\}.
\end{equation*}

Let $N_1,N_2$ be positive integers, where
$N_1=1$ or $N_2=1$ is allowed for one-dimensional
signals. Every dual-quaternion-valued vector or matrix
$\widehat{\mathbf{X}}\in\mathbb{DH}^{N_1\times N_2}$
can be written in the form
\begin{equation}
\widehat{\mathbf{X}}
=
\mathbf{X}_{i1} + \mathbf{X}_{i2}\mathbf{j}
+\varepsilon\!\left(\mathbf{X}_{i3} + \mathbf{X}_{i4}\mathbf{j}\right),
\label{eq:f_ci_decomposition}
\end{equation}
where $\mathbf{X}_{i1}, \mathbf{X}_{i2}, \mathbf{X}_{i3}, \mathbf{X}_{i4}
\in \mathbb C_{\mathbf{i}}^{N_1\times N_2}$.

Similarly, every dual-quaternion-valued vector or matrix
$\widehat{\mathbf{Y}}\in\mathbb{DH}^{N_1\times N_2}$
can be represented as
\begin{equation}
\widehat{\mathbf{Y}} = \mathbf{Y}_{j1} + \mathbf{i} \mathbf{Y}_{j2}
+\varepsilon\!\left(\mathbf{Y}_{j3} + \mathbf{i} \mathbf{Y}_{j4}\right),
\label{eq:g_cj_decomposition}
\end{equation}
where $\mathbf{Y}_{j1}, \mathbf{Y}_{j2}, \mathbf{Y}_{j3}, \mathbf{Y}_{j4}
\in \mathbb C_{\mathbf{j}}^{N_1\times N_2}$.

\begin{remark}[Two decompositions]
\label{rem:two_decompositions}
The decompositions \eqref{eq:f_ci_decomposition} and
\eqref{eq:g_cj_decomposition} differ in which complex subalgebra is
used to absorb the quaternion basis units, and hence in the side on
which the decomposition acts. The first is adapted to the
$\mathbf{i}$-generated left kernel, and the second to the
$\mathbf{j}$-generated right kernel. Because $\mathbf{i}$ and
$\mathbf{j}$ do not commute, no single decomposition is
simultaneously natural for both kernels; the two forms---referred to
below as the $\mathbb C_{\mathbf{i}}$-, $\mathbb C_{\mathbf{j}}$-
decompositions, respectively---are therefore used interchangeably
throughout the paper, according to which kernel acts on the array
first.
\end{remark}

\subsection{Linear Canonical Transform}\label{subsec:lct}

For each dimension $s\in\{1,2\}$, let
\begin{equation*}
A_s=
\begin{bmatrix}
a_s & b_s\\
c_s & d_s
\end{bmatrix},
\qquad a_s d_s-b_s c_s=1,
\end{equation*}
be a real unimodular matrix. In the following, we consider the
nondegenerate case $b_s\neq 0$, $s=1,2$.

Let $N_s$ be a positive integer, and let $\Delta t_s$ and
$\Delta u_s$ denote the uniform sampling steps in the original and
canonical transform domains, respectively. Assume that
\begin{equation*}
\Delta u_s=\frac{2\pi b_s}{N_s\Delta t_s},
\qquad s=1,2.
\end{equation*}
For the unimodular matrix $A_s$, define the discrete canonical phase
by
\begin{equation*}
\Phi_{A_s}(n_s,m_s)
= \frac{a_s}{2b_s}n_s^2\Delta t_s^2
-\frac{2\pi}{N_s}m_s n_s
+\frac{d_s}{2b_s}m_s^2\Delta u_s^2,
\end{equation*}
where $n_s,m_s \in \{0,1,\ldots,N_s-1\}$, $s=1,2$.

Let $\boldsymbol{\mu}$ be a pure unit quaternion satisfying
$\boldsymbol{\mu}^2=-1$. The corresponding discrete LCT kernel is
defined by
\begin{equation}
[\mathbf{K}^{A_s,\boldsymbol{\mu}}]_{n_s,m_s}
=\frac{1}{\sqrt{N_s}}
\exp
\left[
\boldsymbol{\mu}\Phi_{A_s}(n_s,m_s)\right],
\label{eq:DQ-DLCT_kernel}
\end{equation}
so that the associated kernel matrix is
$\mathbf{K}^{A_s,\boldsymbol{\mu}}
=
\frac{1}{\sqrt{N_s}}
\left[e^{\boldsymbol{\mu}\Phi_{A_s}(n_s,m_s)}
\right]_{n_s,m_s=0}^{N_s-1}$.
Since $\Phi_{A_s}(n_s,m_s)$ is real, the quaternion conjugate of the
kernel satisfies
\begin{align*}
[\overline{\mathbf{K}}^{A_s,\boldsymbol{\mu}}]_{n_s,m_s}
&=
\frac{1}{\sqrt{N_s}}
\exp
\left[
-\boldsymbol{\mu}
\Phi_{A_s}(n_s,m_s)
\right]
\nonumber\\
&=
[\mathbf{K}^{A_s,-\boldsymbol{\mu}}]_{n_s,m_s}.
\end{align*}
Accordingly,
$\overline{\mathbf{K}}^{A_s,\boldsymbol{\mu}}
=\mathbf{K}^{A_s,-\boldsymbol{\mu}}$.
For the one-dimensional case, the subscript $s$ is omitted for
simplicity.

\section{Dual Quaternion Discrete Linear Canonical Transform}\label{sec:DQ-DLCT}

\subsection{1-D DQ-DLCTs}
\label{subsec:1d_DQ-DLCT}

Consider a one-dimensional discrete signal
$\widehat{\mathbf{X}}\in\mathbb{DH}^{N\times 1}$, indexed by
$n\in\{0,1,\ldots,N-1\}$. For the one-dimensional case, the
subscript $s$ in the LCT parameters and sampling intervals is
omitted for simplicity.

\begin{definition}[1-D left-sided DQ-DLCT]
The one-dimensional left-sided dual quaternion discrete linear
canonical transform associated with $A$ is defined by
\begin{align}
\mathcal L_{A,\mathbf{i}}^{L}
[\widehat{\mathbf{X}}](m)
&=
\frac{1}{\sqrt{N}}
\sum_{n=0}^{N-1}
\exp
\left[
\mathbf{i}
\left(
\frac{a}{2b}n^2\Delta t^2
-\frac{2\pi}{N}mn
+\frac{d}{2b}m^2\Delta u^2
\right)
\right]
[\widehat{\mathbf{X}}]_n
\nonumber\\
&=
\sum_{n=0}^{N-1}
[\mathbf K^{A,\mathbf{i}}]_{n,m}
[\widehat{\mathbf{X}}]_n
\nonumber\\
&=
\left[
(\mathbf K^{A,\mathbf{i}})^T
\widehat{\mathbf{X}}
\right]_m,
\end{align}
where $m\in\{0,1,\ldots,N-1\}$ and the transpose accounts for the
index ordering $[\mathbf K^{A,\mathbf{i}}]_{n,m}$.
\end{definition}

\begin{definition}[1-D right-sided DQ-DLCT]
The one-dimensional right-sided dual quaternion discrete linear
canonical transform associated with $A$ is defined by
\begin{align}
\mathcal L_{A,\mathbf{i}}^{R}
[\widehat{\mathbf{X}}](m)
&=
\frac{1}{\sqrt{N}}
\sum_{n=0}^{N-1}
[\widehat{\mathbf{X}}]_n
\exp
\left[
\mathbf{i}
\left(
\frac{a}{2b}n^2\Delta t^2
-\frac{2\pi}{N}mn
+\frac{d}{2b}m^2\Delta u^2
\right)
\right]
\nonumber\\
&=
\sum_{n=0}^{N-1}
[\widehat{\mathbf{X}}]_n
[\mathbf K^{A,\mathbf{i}}]_{n,m}
\nonumber\\
&=
\left[
\left(
\widehat{\mathbf{X}}^T
\mathbf K^{A,\mathbf{i}}
\right)^T
\right]_m,
\end{align}
where $m\in\{0,1,\ldots,N-1\}$ and the outer transpose accounts for
the index ordering $[\mathbf K^{A,\mathbf{i}}]_{n,m}$.
\end{definition}

Since multiplication in the dual quaternion algebra is
noncommutative, the left- and right-sided transforms are generally
different:
$\mathcal L_{A,\mathbf{i}}^{L}
[\widehat{\mathbf{X}}]
\neq
\mathcal L_{A,\mathbf{i}}^{R}
[\widehat{\mathbf{X}}]$.
In particular, the relative order of the transform kernel and the
signal samples cannot in general be interchanged. The two transforms
coincide only when the kernel commutes with the signal samples, which
is not the case for a general
$\widehat{\mathbf{X}}\in\mathbb{DH}^{N\times 1}$.

\subsection{2-D DQ-DLCTs}
\label{subsec:2d_DQ-DLCT}

Consider a two-dimensional discrete signal
$\widehat{\mathbf{X}}\in\mathbb{DH}^{N_1\times N_2}$, indexed by
$n_s,m_s\in\{0,1,\ldots,N_s-1\}$, $s=1,2$.

\begin{definition}[2-D TDQ-DLCT]
The two-dimensional two-sided dual quaternion discrete linear
canonical transform (2-D TDQ-DLCT) is defined by
\begin{align}
&\mathcal L_{A_1,\mathbf{i},A_2,\mathbf{j}}^{TS}
[\widehat{\mathbf{X}}](m_1,m_2)
\nonumber\\
&=
\frac{1}{\sqrt{N_1N_2}}
\sum_{n_1=0}^{N_1-1}
\sum_{n_2=0}^{N_2-1}
\exp
\left[
\mathbf{i}
\left(
\frac{a_1}{2b_1}n_1^2\Delta t_1^2
-\frac{2\pi}{N_1}m_1n_1
+\frac{d_1}{2b_1}m_1^2\Delta u_1^2
\right)
\right]
\nonumber\\
&\qquad
[\widehat{\mathbf{X}}]_{n_1,n_2}
\exp
\left[
\mathbf{j}
\left(
\frac{a_2}{2b_2}n_2^2\Delta t_2^2
-\frac{2\pi}{N_2}m_2n_2
+\frac{d_2}{2b_2}m_2^2\Delta u_2^2
\right)
\right]
\nonumber\\
&=
\sum_{n_1=0}^{N_1-1}
\sum_{n_2=0}^{N_2-1}
[\mathbf K^{A_1,\mathbf{i}}]_{n_1,m_1}
[\widehat{\mathbf{X}}]_{n_1,n_2}
[\mathbf K^{A_2,\mathbf{j}}]_{n_2,m_2}
\nonumber\\
&=
\left[
(\mathbf K^{A_1,\mathbf{i}})^T
\widehat{\mathbf{X}}
\mathbf K^{A_2,\mathbf{j}}
\right]_{m_1,m_2},
\label{eq:2d_two_sided_DQ-DLCT}
\end{align}
where the two kernel matrices act on the first and second
dimensions, respectively, and the transpose reflects the index
ordering $[\mathbf K^{A_s,\boldsymbol{\mu}}]_{n_s,m_s}$.
\end{definition}

\begin{definition}[2-D left-sided DQ-DLCT]
The two-dimensional left-sided dual quaternion discrete linear
canonical transform is defined by
\begin{align}
&\mathcal L_{A_1,\mathbf{i},A_2,\mathbf{j}}^{L}
[\widehat{\mathbf{X}}](m_1,m_2)
\nonumber\\
&=
\frac{1}{\sqrt{N_1N_2}}
\sum_{n_1=0}^{N_1-1}
\sum_{n_2=0}^{N_2-1}
\exp
\left[
\mathbf{i}
\left(
\frac{a_1}{2b_1}n_1^2\Delta t_1^2
-\frac{2\pi}{N_1}m_1n_1
+\frac{d_1}{2b_1}m_1^2\Delta u_1^2
\right)
\right]
\nonumber\\
&\qquad
\exp
\left[
\mathbf{j}
\left(
\frac{a_2}{2b_2}n_2^2\Delta t_2^2
-\frac{2\pi}{N_2}m_2n_2
+\frac{d_2}{2b_2}m_2^2\Delta u_2^2
\right)
\right]
[\widehat{\mathbf{X}}]_{n_1,n_2}
\nonumber\\
&=
\sum_{n_1=0}^{N_1-1}
\sum_{n_2=0}^{N_2-1}
[\mathbf K^{A_1,\mathbf{i}}]_{n_1,m_1}
[\mathbf K^{A_2,\mathbf{j}}]_{n_2,m_2}
[\widehat{\mathbf{X}}]_{n_1,n_2}.
\end{align}
\end{definition}
In the left-sided transform, both kernels act from the left of the
signal. Because quaternion multiplication is noncommutative, the
order of the two kernels is fixed:
$\mathbf K^{A_2,\mathbf{j}}$ is applied first, along the second
dimension, followed by $\mathbf K^{A_1,\mathbf{i}}$, along the
first.

\begin{definition}[2-D right-sided DQ-DLCT]
The two-dimensional right-sided dual quaternion discrete linear
canonical transform is defined by
\begin{align}
&\mathcal L_{A_1,\mathbf{i},A_2,\mathbf{j}}^{R}
[\widehat{\mathbf{X}}](m_1,m_2)
\nonumber\\
&=
\frac{1}{\sqrt{N_1N_2}}
\sum_{n_1=0}^{N_1-1}
\sum_{n_2=0}^{N_2-1}
[\widehat{\mathbf{X}}]_{n_1,n_2}
\exp
\left[
\mathbf{i}
\left(
\frac{a_1}{2b_1}n_1^2\Delta t_1^2
-\frac{2\pi}{N_1}m_1n_1
+\frac{d_1}{2b_1}m_1^2\Delta u_1^2
\right)
\right]
\nonumber\\
&\qquad
\exp
\left[
\mathbf{j}
\left(
\frac{a_2}{2b_2}n_2^2\Delta t_2^2
-\frac{2\pi}{N_2}m_2n_2
+\frac{d_2}{2b_2}m_2^2\Delta u_2^2
\right)
\right]
\nonumber\\
&=
\sum_{n_1=0}^{N_1-1}
\sum_{n_2=0}^{N_2-1}
[\widehat{\mathbf{X}}]_{n_1,n_2}
[\mathbf K^{A_1,\mathbf{i}}]_{n_1,m_1}
[\mathbf K^{A_2,\mathbf{j}}]_{n_2,m_2}.
\end{align}
\end{definition}
In the right-sided transform, both kernels act from the right of the
signal. As in the left-sided case, the order of the two kernels is
fixed, with $\mathbf K^{A_1,\mathbf{i}}$ applied first and
$\mathbf K^{A_2,\mathbf{j}}$ second.

Because quaternion multiplication is noncommutative, the order of the
two kernels in the left- and right-sided transforms is fixed. In
general,
$[\mathbf K^{A_1,\mathbf{i}}]_{n_1,m_1}
[\mathbf K^{A_2,\mathbf{j}}]_{n_2,m_2}
\neq
[\mathbf K^{A_2,\mathbf{j}}]_{n_2,m_2}
[\mathbf K^{A_1,\mathbf{i}}]_{n_1,m_1}$.
Hence, exchanging the order of the two directional kernels generally
defines a different transform.

\subsection{Fast Algorithms for DQ-DLCT}
\label{subsec:fast_DQ-DLCT}

For each dimension $s\in\{1,2\}$ and a pure unit quaternion
$\boldsymbol{\mu}$ satisfying $\boldsymbol{\mu}^2=-1$, define the
normalized quaternion Fourier matrix by
\begin{equation*}
[\mathbf F_{N_s}^{\boldsymbol{\mu}}]_{n_s,m_s}
=
\frac{1}{\sqrt{N_s}}
\exp
\left(
-\boldsymbol{\mu}
\frac{2\pi}{N_s}
n_sm_s
\right),
\qquad
n_s,m_s\in\{0,1,\ldots,N_s-1\}.
\end{equation*}

We use $\mathbf P_s^{\boldsymbol{\mu}}$ and
$\mathbf Q_s^{\boldsymbol{\mu}}$ to denote the input- and
output-domain chirp matrices, respectively:
\begin{align}
\mathbf P_s^{\boldsymbol{\mu}}
&=
\operatorname{diag}
\left\{
\exp
\left(
\boldsymbol{\mu}
\frac{a_s}{2b_s}
n_s^2\Delta t_s^2
\right)
\right\}_{n_s=0}^{N_s-1},
\\
\mathbf Q_s^{\boldsymbol{\mu}}
&=
\operatorname{diag}
\left\{
\exp
\left(
\boldsymbol{\mu}
\frac{d_s}{2b_s}
m_s^2\Delta u_s^2
\right)
\right\}_{m_s=0}^{N_s-1}.
\end{align}

The following factorization of the canonical kernel is the basis
for the fast algorithms developed in this section.

\begin{theorem}[Chirp--Fourier--chirp factorization]
\label{thm:kernel_factorization}
For each $s\in\{1,2\}$ and each pure unit quaternion
$\boldsymbol{\mu}$ with $\boldsymbol{\mu}^2=-1$, the discrete
canonical kernel admits the factorization
\begin{equation}
\mathbf K^{A_s,\boldsymbol{\mu}}
=
\mathbf P_s^{\boldsymbol{\mu}}
\mathbf F_{N_s}^{\boldsymbol{\mu}}
\mathbf Q_s^{\boldsymbol{\mu}}.
\label{eq:DQ-DLCT_kernel_factorization}
\end{equation}
\end{theorem}

\begin{proof}
The $(n_s,m_s)$ entry of the right-hand side equals
\begin{align*}
&
\left[
\mathbf P_s^{\boldsymbol{\mu}}
\mathbf F_{N_s}^{\boldsymbol{\mu}}
\mathbf Q_s^{\boldsymbol{\mu}}
\right]_{n_s,m_s}
\nonumber\\
&=
\frac{1}{\sqrt{N_s}}
\exp
\left(
\boldsymbol{\mu}
\frac{a_s}{2b_s}
n_s^2\Delta t_s^2
\right)
\exp
\left(
-\boldsymbol{\mu}
\frac{2\pi}{N_s}
n_sm_s
\right)
\exp
\left(
\boldsymbol{\mu}
\frac{d_s}{2b_s}
m_s^2\Delta u_s^2
\right)
\nonumber\\
&=
\frac{1}{\sqrt{N_s}}
\exp
\left[
\boldsymbol{\mu}
\left(
\frac{a_s}{2b_s}n_s^2\Delta t_s^2
-\frac{2\pi}{N_s}n_sm_s
+\frac{d_s}{2b_s}m_s^2\Delta u_s^2
\right)
\right]
\nonumber\\
&=
[\mathbf K^{A_s,\boldsymbol{\mu}}]_{n_s,m_s}.
\end{align*}
Here the scalar exponential factors are generated by the same pure
unit quaternion $\boldsymbol{\mu}$ and therefore commute. No
commutativity of the matrix factors themselves is required.
\end{proof}

\begin{remark}[Computational cost]
\label{rem:fast_cost}
Since $\mathbf P_s^{\boldsymbol{\mu}}$ and
$\mathbf Q_s^{\boldsymbol{\mu}}$ are diagonal and
$\mathbf F_{N_s}^{\boldsymbol{\mu}}$ can be computed using complex
FFTs, the factorization \eqref{eq:DQ-DLCT_kernel_factorization}
reduces the cost of applying $\mathbf K^{A_s,\boldsymbol{\mu}}$
from $\mathcal O(N_s^2)$ to $\mathcal O(N_s\log N_s)$ per
dimension.
\end{remark}

For the one-dimensional case, the subscript $s$ is omitted, so that
\begin{equation*}
\mathbf K^{A,\boldsymbol{\mu}}
=
\mathbf P^{\boldsymbol{\mu}}
\mathbf F_N^{\boldsymbol{\mu}}
\mathbf Q^{\boldsymbol{\mu}}.
\end{equation*}

\subsubsection{Fast 1-D DQ-DLCT}
\label{subsubsec:fast_1d_DQ-DLCT}

Consider first the 1-D left-sided DQ-DLCT. Using the
$\mathbb C_{\mathbf{i}}$-decomposition introduced in
Section~\ref{subsec:dual_quaternion}, write
$\widehat{\mathbf{X}}
=
\mathbf{X}_{i1}
+\mathbf{X}_{i2}\mathbf{j}
+\varepsilon
\left(\mathbf{X}_{i3}+\mathbf{X}_{i4}\mathbf{j}\right)$,
with $\mathbf{X}_{i1},\mathbf{X}_{i2},\mathbf{X}_{i3},\mathbf{X}_{i4}
\in
\mathbb C_{\mathbf{i}}^{N\times1}$.

Since the normalized Fourier matrix is symmetric,
$\left(\mathbf F_N^{\mathbf{i}}\right)^T
=
\mathbf F_N^{\mathbf{i}}$,
and both chirp matrices are diagonal, it follows that
$\left(
\mathbf K^{A,\mathbf{i}}
\right)^T
=
\mathbf Q^{\mathbf{i}}
\mathbf F_N^{\mathbf{i}}
\mathbf P^{\mathbf{i}}$.

Therefore,
\begin{align*}
\mathcal L_{A,\mathbf{i}}^{L}
[\widehat{\mathbf{X}}]
&=
\mathbf Q^{\mathbf{i}}
\mathbf F_N^{\mathbf{i}}
\mathbf P^{\mathbf{i}}
\widehat{\mathbf{X}}
\nonumber\\
&=
\mathbf Q^{\mathbf{i}}
\Bigg[
\mathbf F_N^{\mathbf{i}}
\left(
\mathbf P^{\mathbf{i}}
\mathbf{X}_{i1}
\right)
+
\mathbf F_N^{\mathbf{i}}
\left(
\mathbf P^{\mathbf{i}}
\mathbf{X}_{i2}
\right)\mathbf{j}
+
\varepsilon
\left(
\mathbf F_N^{\mathbf{i}}
\left(
\mathbf P^{\mathbf{i}}
\mathbf{X}_{i3}
\right)
+
\mathbf F_N^{\mathbf{i}}
\left(
\mathbf P^{\mathbf{i}}
\mathbf{X}_{i4}
\right)\mathbf{j}
\right)
\Bigg].
\end{align*}

Thus, the 1-D left-sided DQ-DLCT can be evaluated through
the following three operations:
\begin{enumerate}
    \item multiply each of the four $\mathbb C_{\mathbf{i}}$-valued
    components by the input chirp $\mathbf P^{\mathbf{i}}$;
    \item apply the normalized complex Fourier transform generated by
    $\mathbf{i}$ to the four chirp-modulated components;
    \item recombine the transformed components and multiply the result
    by the output chirp $\mathbf Q^{\mathbf{i}}$.
\end{enumerate}

The chirp multiplications and component recombination require
$\mathcal O(N)$ operations. The four Fourier transforms share the
same transform matrix $\mathbf F_N^{\mathbf{i}}$ and can be computed
by four complex FFTs, or batched into a single call, at a cost of
$\mathcal O(N\log N)$ operations. By contrast, direct evaluation of
the defining summation requires $\mathcal O(N^2)$ operations.

For the 1-D right-sided DQ-DLCT, the components containing a right
factor $\mathbf{j}$ require special treatment because they do not
commute with the $\mathbf{i}$-generated kernel. For every
$z\in\mathbb C_{\mathbf{i}}$,
$
\mathbf{j} z
=
\overline{z}\,\mathbf{j}.
$
In particular,
$
\mathbf{j}
e^{\mathbf{i}\theta}
=
e^{-\mathbf{i}\theta}
\mathbf{j}.
$

Using
$
\overline{\mathbf P^{\mathbf{i}}}
=
\mathbf P^{-\mathbf{i}},
\overline{\mathbf F_N^{\mathbf{i}}}
=
\mathbf F_N^{-\mathbf{i}},
\overline{\mathbf Q^{\mathbf{i}}}
=
\mathbf Q^{-\mathbf{i}},$
the 1-D right-sided DQ-DLCT can be written in column-vector form as
\begin{align*}
\mathcal L_{A,\mathbf{i}}^{R}
[\widehat{\mathbf{X}}]
&=
\mathbf Q^{\mathbf{i}}
\mathbf F_N^{\mathbf{i}}
\mathbf P^{\mathbf{i}}
\mathbf{X}_{i1}
+
\mathbf Q^{-\mathbf{i}}
\mathbf F_N^{-\mathbf{i}}
\mathbf P^{-\mathbf{i}}
\mathbf{X}_{i2}\mathbf{j}
\nonumber\\
&\quad
+
\varepsilon
\Big(
\mathbf Q^{\mathbf{i}}
\mathbf F_N^{\mathbf{i}}
\mathbf P^{\mathbf{i}}
\mathbf{X}_{i3}
+
\mathbf Q^{-\mathbf{i}}
\mathbf F_N^{-\mathbf{i}}
\mathbf P^{-\mathbf{i}}
\mathbf{X}_{i4}\mathbf{j}
\Big).
\end{align*}

Hence, the two components without a right factor $\mathbf{j}$ are
evaluated with the $\mathbf{i}$-generated Fourier and chirp factors,
whereas the two components containing $\mathbf{j}$ are evaluated with
their conjugated counterparts. Since
$\mathbf F_N^{-\mathbf{i}} =
\overline{\mathbf F_N^{\mathbf{i}}}$, both transforms can be
implemented using standard complex FFT routines, and the four
Fourier applications can be grouped into two pairs sharing the same
transform matrix. Consequently, the right-sided DQ-DLCT also has
complexity $\mathcal O(N\log N)$.

The analogous fast algorithm for the two-sided and multidimensional
cases follows the same chirp--Fourier--chirp pattern and is developed
in Section~\ref{subsubsec:fast_2d_tDQ-DLCT}.
\subsubsection{Fast 2-D TDQ-DLCT}
\label{subsubsec:fast_2d_tDQ-DLCT}
We next consider $\widehat{\mathbf{X}}\in\mathbb{DH}^{N_1\times N_2}$.
Applying the kernel factorization
\eqref{eq:DQ-DLCT_kernel_factorization} to the two dimensions gives
\begin{equation}
\mathbf K^{A_1,\mathbf{i}}
=
\mathbf P_1^{\mathbf{i}}
\mathbf F_{N_1}^{\mathbf{i}}
\mathbf Q_1^{\mathbf{i}},
\quad
\mathbf K^{A_2,\mathbf{j}}
=
\mathbf P_2^{\mathbf{j}}
\mathbf F_{N_2}^{\mathbf{j}}
\mathbf Q_2^{\mathbf{j}}.
\label{eq:2d_kernel_factorization}
\end{equation}

Since $\mathbf F_{N_1}^{\mathbf{i}}$ is symmetric and
$\mathbf P_1^{\mathbf{i}}$ and $\mathbf Q_1^{\mathbf{i}}$ are
diagonal,
\begin{equation}
\left(
\mathbf K^{A_1,\mathbf{i}}
\right)^T
=
\mathbf Q_1^{\mathbf{i}}
\mathbf F_{N_1}^{\mathbf{i}}
\mathbf P_1^{\mathbf{i}}.
\label{eq:left_kernel_fft_factorization}
\end{equation}
Substituting \eqref{eq:2d_kernel_factorization}
and \eqref{eq:left_kernel_fft_factorization} into the matrix
representation of the TDQ-DLCT yields
\begin{align}
\mathcal L_{A_1,\mathbf{i},A_2,\mathbf{j}}^{TS}
[\widehat{\mathbf{X}}]
&=
\left(
\mathbf K^{A_1,\mathbf{i}}
\right)^T
\widehat{\mathbf{X}}
\mathbf K^{A_2,\mathbf{j}}
\nonumber\\
&=
\mathbf Q_1^{\mathbf{i}}
\mathbf F_{N_1}^{\mathbf{i}}
\mathbf P_1^{\mathbf{i}}
\widehat{\mathbf{X}}
\mathbf P_2^{\mathbf{j}}
\mathbf F_{N_2}^{\mathbf{j}}
\mathbf Q_2^{\mathbf{j}}.
\label{eq:2d_DQ-DLCT_factorization}
\end{align}

Define the input-chirp-modulated signal
$
\widehat{\mathbf{Y}}
:=
\mathbf P_1^{\mathbf{i}}
\widehat{\mathbf{X}}
\mathbf P_2^{\mathbf{j}}.
$
Then
\begin{equation}
\mathcal L_{A_1,\mathbf{i},A_2,\mathbf{j}}^{TS}
[\widehat{\mathbf{X}}]
=
\mathbf Q_1^{\mathbf{i}}
\left(
\mathbf F_{N_1}^{\mathbf{i}}
\widehat{\mathbf{Y}}
\mathbf F_{N_2}^{\mathbf{j}}
\right)
\mathbf Q_2^{\mathbf{j}}.
\label{eq:DQ-DLCT_as_dqdft}
\end{equation}

Equation~\eqref{eq:DQ-DLCT_as_dqdft} shows that the TDQ-DLCT consists
of three stages: input chirp modulation, a normalized two-sided dual
quaternion discrete Fourier transform (TDQ-DFT), and output chirp
modulation. Once the chirp factors have been applied, the dominant
computation is the central two-sided Fourier transform, which is
evaluated by the fast TDQ-DFT procedure developed in
Section~\ref{subsubsec:fast_2d_tDQ-DLCT}.

To evaluate this Fourier stage efficiently, decompose the
chirp-modulated signal as
\begin{equation*}
\widehat{\mathbf{Y}}
=
\mathbf{Y}_{i1}
+
\mathbf{Y}_{i2}\mathbf{j}
+
\varepsilon
\left(
\mathbf{Y}_{i3}
+
\mathbf{Y}_{i4}\mathbf{j}
\right),
\qquad
\mathbf{Y}_{i1},\mathbf{Y}_{i2},\mathbf{Y}_{i3},\mathbf{Y}_{i4}
\in
\mathbb C_{\mathbf{i}}^{N_1\times N_2}.
\end{equation*}

Since every entry of $\mathbf F_{N_2}^{\mathbf{j}}$ belongs to
$\mathbb C_{\mathbf{j}}$, it commutes with $\mathbf{j}$. Hence,
\begin{align*}
\mathbf F_{N_1}^{\mathbf{i}}
\widehat{\mathbf{Y}}
\mathbf F_{N_2}^{\mathbf{j}}
&=
\mathbf F_{N_1}^{\mathbf{i}}
\mathbf{Y}_{i1}
\mathbf F_{N_2}^{\mathbf{j}}
+
\mathbf F_{N_1}^{\mathbf{i}}
\mathbf{Y}_{i2}
\mathbf F_{N_2}^{\mathbf{j}}
\mathbf{j}
\nonumber\\
&\quad
+
\varepsilon
\Big(
\mathbf F_{N_1}^{\mathbf{i}}
\mathbf{Y}_{i3}
\mathbf F_{N_2}^{\mathbf{j}}
+
\mathbf F_{N_1}^{\mathbf{i}}
\mathbf{Y}_{i4}
\mathbf F_{N_2}^{\mathbf{j}}
\mathbf{j}
\Big).
\end{align*}
The above can be evaluated by the fast TDQ-DFT procedure using a
fixed number of ordinary complex FFTs. The factorization
\eqref{eq:DQ-DLCT_as_dqdft} then leads directly to the computational
procedure given in Algorithm~\ref{alg:fast_tDQ-DLCT}.

\begin{algorithm}[t]
\caption{Fast two-dimensional TDQ-DLCT}
\label{alg:fast_tDQ-DLCT}
\begin{algorithmic}[1]
\Require $\widehat{\mathbf{X}}$, $A_1$, $A_2$, $\Delta t_1$, $\Delta t_2$
\Ensure $\widehat{\mathbf{Y}}=\mathcal L_{A_1,\mathbf{i},A_2,\mathbf{j}}^{TS}[\widehat{\mathbf{X}}]$
\State Construct $\mathbf P_1^{\mathbf{i}}$, $\mathbf Q_1^{\mathbf{i}}$,
$\mathbf P_2^{\mathbf{j}}$, and $\mathbf Q_2^{\mathbf{j}}$.
\State $\widehat{\mathbf{Y}}\gets
\mathbf P_1^{\mathbf{i}}\widehat{\mathbf{X}}\mathbf P_2^{\mathbf{j}}$.
\State $\widehat{\mathbf{H}}\gets
\mathbf F_{N_1}^{\mathbf{i}}\widehat{\mathbf{Y}}\mathbf F_{N_2}^{\mathbf{j}}$
using the fast two-sided DQ-DFT procedure of
Section~\ref{subsubsec:fast_2d_tDQ-DLCT}.
\State $\widehat{\mathbf{Y}}\gets
\mathbf Q_1^{\mathbf{i}}\widehat{\mathbf{H}}\mathbf Q_2^{\mathbf{j}}$.
\State \Return $\widehat{\mathbf{Y}}$.
\end{algorithmic}
\end{algorithm}

The TDQ-DLCT is therefore implemented as input chirp modulation, a
fast TDQ-DFT, and output chirp modulation.

Since the chirp matrices are diagonal, constructing them costs
$\mathcal O(N_1+N_2)$ operations, and the chirp operations
themselves require $\mathcal O(N_1N_2)$ operations. The dominant
cost lies in the Fourier stage: with a row--column FFT
implementation, in which $\mathbf F_{N_1}^{\mathbf{i}}$ is applied
along the first dimension and $\mathbf F_{N_2}^{\mathbf{j}}$ along
the second, this stage has complexity
$
\mathcal O\!\left(
N_1N_2(\log N_1+\log N_2)
\right).
$
By contrast, direct evaluation requires
$\mathcal O(N_1^2N_2^2)$ operations.

The inverse TDQ-DLCT admits the same factorized structure. Let
$
\widehat{\mathbf{Y}}
=
\mathcal L_{A_1,\mathbf{i},A_2,\mathbf{j}}^{TS}
[\widehat{\mathbf{X}}].
$
By \eqref{eq:2d_DQ-DLCT_factorization}, the inverse transform is
\begin{align}
\widehat{\mathbf{X}}
&=
\left(\mathbf P_1^{\mathbf{i}}\right)^{-1}
\left(\mathbf F_{N_1}^{\mathbf{i}}\right)^{-1}
\left(\mathbf Q_1^{\mathbf{i}}\right)^{-1}
\widehat{\mathbf{Y}}
\left(\mathbf Q_2^{\mathbf{j}}\right)^{-1}
\left(\mathbf F_{N_2}^{\mathbf{j}}\right)^{-1}
\left(\mathbf P_2^{\mathbf{j}}\right)^{-1}.
\end{align}
Since the chirp matrices have unit-modulus diagonal entries,
\begin{equation*}
\left(\mathbf P_s^{\boldsymbol{\mu}}\right)^{-1}
=
\mathbf P_s^{-\boldsymbol{\mu}},
\qquad
\left(\mathbf Q_s^{\boldsymbol{\mu}}\right)^{-1}
=
\mathbf Q_s^{-\boldsymbol{\mu}}.
\end{equation*}
Thus the inverse chirp factors are obtained by reversing the phase
signs, while the inverse Fourier stage is computed by applying
inverse FFTs along each dimension. The inverse transform can
therefore be implemented as described in
Algorithm~\ref{alg:fast_inverse_tDQ-DLCT} and inherits the same
asymptotic complexity,
$\mathcal{O}\!\left(N_1N_2(\log N_1+\log N_2)\right)$.

\begin{algorithm}[t]
\caption{Fast inverse 2-D TDQ-DLCT}
\label{alg:fast_inverse_tDQ-DLCT}
\begin{algorithmic}[1]
\Require $\widehat{\mathbf{Y}}$, $A_1$, $A_2$, $\Delta t_1$, $\Delta t_2$
\Ensure $\widehat{\mathbf{X}}$
\State $\widehat{\mathbf{U}}\gets
(\mathbf Q_1^{\mathbf{i}})^{-1}\widehat{\mathbf{Y}}
(\mathbf Q_2^{\mathbf{j}})^{-1}$.
\State $\widehat{\mathbf{V}}\gets
(\mathbf F_{N_1}^{\mathbf{i}})^{-1}\widehat{\mathbf{U}}
(\mathbf F_{N_2}^{\mathbf{j}})^{-1}$ using the inverse fast two-sided DQ-DFT.
\State $\widehat{\mathbf{X}}\gets
(\mathbf P_1^{\mathbf{i}})^{-1}\widehat{\mathbf{V}}
(\mathbf P_2^{\mathbf{j}})^{-1}$.
\State \Return $\widehat{\mathbf{X}}$.
\end{algorithmic}
\end{algorithm}


\section{Numerical Properties of the TDQ-DLCT}
\label{sec:numerical_properties}

The factorized representation of the discrete canonical kernel also
makes its numerical structure explicit. For a quaternion matrix
$\mathbf M$, let $\mathbf M^{H}=\overline{\mathbf M}^{\,T}$ denote
the conjugate transpose. For brevity, we write $\mathcal L^{TS}$ for
$\mathcal L_{A_1,\mathbf{i},A_2,\mathbf{j}}^{TS}$ when the parameters
are clear from context.

\begin{theorem}[Unitarity of the discrete canonical kernel]
\label{thm:kernel_unitarity}
For each $s\in\{1,2\}$ and each pure unit quaternion
$\boldsymbol{\mu}$ with $\boldsymbol{\mu}^2=-1$, the kernel matrix
$\mathbf K^{A_s,\boldsymbol{\mu}}$ defined in
\eqref{eq:DQ-DLCT_kernel} is unitary:
\begin{equation}
\left(\mathbf K^{A_s,\boldsymbol{\mu}}\right)^H
\mathbf K^{A_s,\boldsymbol{\mu}}
=
\mathbf K^{A_s,\boldsymbol{\mu}}
\left(\mathbf K^{A_s,\boldsymbol{\mu}}\right)^H
=
\mathbf I_{N_s}.
\label{eq:kernel_unitarity}
\end{equation}
\end{theorem}

\begin{proof}
By Theorem~\ref{thm:kernel_factorization},
$\mathbf K^{A_s,\boldsymbol{\mu}}
=
\mathbf P_s^{\boldsymbol{\mu}}
\mathbf F_{N_s}^{\boldsymbol{\mu}}
\mathbf Q_s^{\boldsymbol{\mu}}$.
All entries of these three matrices lie in the commutative complex
subalgebra generated by $\boldsymbol{\mu}$, so the ordinary rules of
complex matrix algebra apply to their products.

The chirp matrices are unitary. Since their diagonal entries have unit
modulus,
\[
\left(\mathbf P_s^{\boldsymbol{\mu}}\right)^H
\mathbf P_s^{\boldsymbol{\mu}}
=
\mathbf I_{N_s},
\qquad
\left(\mathbf Q_s^{\boldsymbol{\mu}}\right)^H
\mathbf Q_s^{\boldsymbol{\mu}}
=
\mathbf I_{N_s}.
\]
The normalized Fourier matrix is also unitary,
\[
\left(\mathbf F_{N_s}^{\boldsymbol{\mu}}\right)^H
\mathbf F_{N_s}^{\boldsymbol{\mu}}
=
\mathbf I_{N_s}.
\]
Substituting the factorization and using the commutativity of the
subalgebra,
\[
\left(\mathbf K^{A_s,\boldsymbol{\mu}}\right)^H
\mathbf K^{A_s,\boldsymbol{\mu}}
=
\left(\mathbf Q_s^{\boldsymbol{\mu}}\right)^H
\left(\mathbf F_{N_s}^{\boldsymbol{\mu}}\right)^H
\left(\mathbf P_s^{\boldsymbol{\mu}}\right)^H
\mathbf P_s^{\boldsymbol{\mu}}
\mathbf F_{N_s}^{\boldsymbol{\mu}}
\mathbf Q_s^{\boldsymbol{\mu}}
=
\mathbf I_{N_s}.
\]
The second identity in \eqref{eq:kernel_unitarity} follows by an
analogous argument from the reverse product.
\end{proof}

For a quaternion-valued array $\mathbf{X}$, define
\begin{equation*}
\|\mathbf{X}\|_F^2
=
\sum_{n_1=0}^{N_1-1}
\sum_{n_2=0}^{N_2-1}
|\mathbf{X}(n_1,n_2)|^2.
\end{equation*}
For a dual quaternion array
$\widehat{\mathbf{X}}=\mathbf{X}_r+\varepsilon \mathbf{X}_d$, with
$\mathbf{X}_r,\mathbf{X}_d\in\mathbb H^{N_1\times N_2}$, define the coefficient energy
norm by
\begin{equation*}
\|\widehat{\mathbf{X}}\|_{\mathcal E}^2
=
\|\mathbf{X}_r\|_F^2+\|\mathbf{X}_d\|_F^2.
\end{equation*}
This norm measures the total energy stored in the primal and dual
coefficients of the dual quaternion array.

\begin{theorem}[Energy preservation and conditioning]
\label{thm:energy_preservation}
Let $\widehat{\mathbf{X}}\in\mathbb{DH}^{N_1\times N_2}$. Then
\begin{equation}
\left\|
\mathcal L_{A_1,\mathbf{i},A_2,\mathbf{j}}^{TS}[\widehat{\mathbf{X}}]
\right\|_{\mathcal E}
=
\|\widehat{\mathbf{X}}\|_{\mathcal E}.
\label{eq:tDQ-DLCT_energy_preservation}
\end{equation}
Consequently, the exact TDQ-DLCT, that is, the transform evaluated in
exact arithmetic, has condition number one with respect to
$\|\cdot\|_{\mathcal E}$:
\begin{equation}
\kappa_{\mathcal E}
\left(\mathcal L_{A_1,\mathbf{i},A_2,\mathbf{j}}^{TS}\right)
=1.
\label{eq:tDQ-DLCT_condition_number}
\end{equation}
\end{theorem}

\begin{proof}
Because $\varepsilon$ commutes with quaternion-valued quantities, the
TDQ-DLCT acts separately on the primal and dual parts:
\[
\mathcal L^{TS}[\widehat{\mathbf{X}}]
=
\mathcal L^{TS}[\mathbf{X}_r]
+\varepsilon\mathcal L^{TS}[\mathbf{X}_d].
\]
For either quaternion component $\mathbf{X}\in\mathbb H^{N_1\times N_2}$,
\[
\mathcal L^{TS}[\mathbf{X}]
=
\left(\mathbf K^{A_1,\mathbf{i}}\right)^T
\mathbf{X}
\mathbf K^{A_2,\mathbf{j}}.
\]
By Theorem~\ref{thm:kernel_unitarity},
$\mathbf K^{A_2,\mathbf{j}}$ is unitary. Moreover, by
\eqref{eq:left_kernel_fft_factorization},
$\left(\mathbf K^{A_1,\mathbf{i}}\right)^T$ is a product of unitary
matrices in $\mathbb C_{\mathbf{i}}$ and is therefore unitary. Since
the quaternion Frobenius norm is invariant under left and right
multiplication by unitary quaternion matrices,
\[
\left\|
\left(\mathbf K^{A_1,\mathbf{i}}\right)^T
\mathbf{X}
\mathbf K^{A_2,\mathbf{j}}
\right\|_F
=
\|\mathbf{X}\|_F.
\]
Applying this identity to $\mathbf{X}_r$ and $\mathbf{X}_d$ and summing the two squared
norms gives \eqref{eq:tDQ-DLCT_energy_preservation}. Consequently, the
induced operator norms of the exact TDQ-DLCT and of its inverse both
equal one, which yields \eqref{eq:tDQ-DLCT_condition_number}.
\end{proof}

\begin{remark}[Why the dual unit does not break unitarity]
\label{rem:dual_unit_unitarity}
Although $\varepsilon$ is nilpotent, it is central and satisfies
$\varepsilon^2=0$. Centrality is what matters here: it allows the
TDQ-DLCT to pass through $\varepsilon$ and act on the primal and dual
components separately, while nilpotency ensures that no cross terms
such as $\varepsilon^2$ arise. Consequently, the total energy is the
sum of the two component energies, which is why the condition number
in \eqref{eq:tDQ-DLCT_condition_number} remains exactly one despite the
presence of the dual unit; no cross terms between primal and dual
parts appear, so no cancellation can occur.
\end{remark}

\begin{corollary}[Exact perturbation stability]
\label{cor:perturbation_stability}
For any perturbation
$\delta\widehat{\mathbf{X}}\in\mathbb{DH}^{N_1\times N_2}$,
\begin{equation}
\left\|
\mathcal L^{TS}[\widehat{\mathbf{X}}+\delta\widehat{\mathbf{X}}]
-
\mathcal L^{TS}[\widehat{\mathbf{X}}]
\right\|_{\mathcal E}
=
\|\delta\widehat{\mathbf{X}}\|_{\mathcal E}.
\label{eq:absolute_perturbation_stability}
\end{equation}
If $\widehat{\mathbf{X}}\neq0$, then
\begin{equation}
\frac{
\left\|
\mathcal L^{TS}[\widehat{\mathbf{X}}+\delta\widehat{\mathbf{X}}]
-
\mathcal L^{TS}[\widehat{\mathbf{X}}]
\right\|_{\mathcal E}
}{
\left\|\mathcal L^{TS}[\widehat{\mathbf{X}}]\right\|_{\mathcal E}
}
=
\frac{\|\delta\widehat{\mathbf{X}}\|_{\mathcal E}}{\|\widehat{\mathbf{X}}\|_{\mathcal E}}.
\label{eq:relative_perturbation_stability}
\end{equation}
In particular, the relative error in the transformed signal equals
the relative error in the input, with no amplification.
\end{corollary}

\begin{proof}
The TDQ-DLCT is linear over addition, so the left-hand side of
\eqref{eq:absolute_perturbation_stability} equals
$\|\mathcal L^{TS}[\delta\widehat{\mathbf{X}}]\|_{\mathcal E}$.
Equation~\eqref{eq:tDQ-DLCT_energy_preservation} then gives
\eqref{eq:absolute_perturbation_stability}, and division by
$\|\mathcal L^{TS}[\widehat{\mathbf{X}}]\|_{\mathcal E}
=\|\widehat{\mathbf{X}}\|_{\mathcal E}$ gives
\eqref{eq:relative_perturbation_stability}.
\end{proof}

\section{Convolution Theorems}\label{sec:convolution}

The classical convolution theorem relates spatial-domain convolution
to multiplication in the transform domain. For dual quaternion
signals, this relation is more involved because multiplication is
noncommutative and the TDQ-DLCT contains left- and right-sided
quadratic phase factors. We therefore introduce spatial and frequency
convolutions compatible with the TDQ-DLCT and derive their
corresponding transform-domain representations.

For two dual-quaternion-valued arrays $\widehat{\mathbf{X}}$ and $\widehat{\mathbf{Y}}$ of the same size,
define their ordered Hadamard product by
\begin{equation}
[\widehat{\mathbf{X}}\circ \widehat{\mathbf{Y}}](m_1,m_2)
=
\widehat{\mathbf{X}}(m_1,m_2)\widehat{\mathbf{Y}}(m_1,m_2).
\label{eq:ordered_hadamard_product}
\end{equation}
Since dual quaternion multiplication is noncommutative, in general
$\widehat{\mathbf{X}}\circ \widehat{\mathbf{Y}}\neq \widehat{\mathbf{Y}}\circ \widehat{\mathbf{X}}$. Hence the multiplication order in
\eqref{eq:ordered_hadamard_product} is fixed throughout this section.
When $\widehat{\mathbf{X}}$ and $\widehat{\mathbf{Y}}$ are scalar-valued, or more generally when they
commute, \eqref{eq:ordered_hadamard_product} reduces to the usual
Hadamard product.

\begin{remark}[Origin of the four component operators]
\label{rem:four_operators}
The four spatial operators
$\mathcal C_{00}, \mathcal C_{01}, \mathcal C_{10}, \mathcal C_{11}$
(and their frequency-domain counterparts
$\mathcal D_{00}, \mathcal D_{01}, \mathcal D_{10}, \mathcal D_{11}$,
defined in Section~\ref{subsec:frequency_convolution})
correspond to the four distinct ways in which the
$\mathbf{i}$-generated and $\mathbf{j}$-generated canonical phases,
as they appear in the kernel factorization, can interact. Because
$\mathbf{i}$ and $\mathbf{j}$ do not commute, the product of the two
phase factors does not simply combine into a single quadratic phase;
instead, the commutation identity
\eqref{eq:kernel_commutation_identity} splits each product into two
terms. The four operators correspond to the four quaternion-component
pairings induced by the $\mathbb C_{\mathbf{i}}$ and
$\mathbb C_{\mathbf{j}}$ decompositions, and they account for the
noncommutative interaction of the left and right canonical kernels.
Throughout, the primal and dual contributions are assembled
separately, since $\varepsilon^2=0$ prevents any mixing between them.
\end{remark}

\subsection{Spatial Convolution}

The following convolution is designed so that the TDQ-DLCT converts
it into a structured sum of ordered Hadamard products with
noncommutative correction terms.

\begin{definition}[Spatial convolution]
Let
$\widehat{\mathbf{X}},\widehat{\mathbf{Y}}
\in\mathbb{DH}^{N_1\times N_2}$.
For an integer $\ell$, let $[\ell]_N$ denote its representative in
$\{0,1,\ldots,N-1\}$, and define
\[
y_1=[n_1-z_1]_{N_1},
\qquad
y_2=[n_2-z_2]_{N_2}.
\]
The TDQ-DLCT-associated spatial convolution is defined by
\begin{align}
(\widehat{\mathbf{X}}*\widehat{\mathbf{Y}})(n_1,n_2)
&=
\sum_{z_1=0}^{N_1-1}
\sum_{z_2=0}^{N_2-1}
\exp\!\left[
\mathbf{i}\frac{a_1}{2b_1}
\left(z_1^2+y_1^2-n_1^2\right)
\Delta t_1^2
\right]
\nonumber\\
&\quad
[\widehat{\mathbf{X}}]_{z_1,z_2}
[\widehat{\mathbf{Y}}]_{y_1,y_2}
\exp\!\left[
\mathbf{j}\frac{a_2}{2b_2}
\left(z_2^2+y_2^2-n_2^2\right)
\Delta t_2^2
\right].
\label{eq:spatial_convolution}
\end{align}
\end{definition}

All spatial indices are interpreted modulo $N_1$ and $N_2$, and the
multiplication order in \eqref{eq:spatial_convolution} is fixed.
When no wrap-around occurs, $y_r=n_r-z_r$ and, for $r=1,2$,
\[
\frac{a_r}{2b_r}
\left(z_r^2+y_r^2-n_r^2\right)\Delta t_r^2
=
\frac{a_r}{b_r}
z_r(z_r-n_r)\Delta t_r^2,
\]
so the cyclic definition locally reduces to the usual chirp form.

For
$\mathbf{X}\in\mathbb C_{\mathbf{i}}^{N_1\times N_2}$ and
$\mathbf{Y}\in\mathbb C_{\mathbf{j}}^{N_1\times N_2}$, define
\begin{align*}
\mathcal C_{00}(\mathbf{X},\mathbf{Y})
&=
(\mathbf K^{A_1,\mathbf{i}})^T
\mathbf{X}
\overline{\mathbf K}^{A_2,\mathbf{j}}
\circ
(\mathbf K^{A_1,\mathbf{i}})^T
\mathbf{Y}
\mathbf K^{A_2,\mathbf{j}}
\nonumber\\
&\quad
+
2(\mathbf K^{A_1,\mathbf{i}})^T
\mathbf{X}\mathbf S
\circ
\mathbf C
\mathbf{Y}
\mathbf K^{A_2,\mathbf{j}}
\mathbf{j},
\\
\mathcal C_{01}(\mathbf{X},\mathbf{Y})
&=
\mathbf{i}
(\mathbf K^{A_1,\mathbf{i}})^T
\mathbf{X}
\overline{\mathbf K}^{A_2,\mathbf{j}}
\circ
(\mathbf K^{A_1,\mathbf{i}})^T
\mathbf{Y}
\mathbf K^{A_2,\mathbf{j}}
\nonumber\\
&\quad
+
2\mathbf{i}
(\mathbf K^{A_1,\mathbf{i}})^T
\mathbf{X}\mathbf S
\circ
\mathbf C
\mathbf{Y}
\mathbf K^{A_2,\mathbf{j}}
\mathbf{j},
\\
\mathcal C_{10}(\mathbf{X},\mathbf{Y})
&=
(\mathbf K^{A_1,\mathbf{i}})^T
\mathbf{X}
\overline{\mathbf K}^{A_2,\mathbf{j}}
\circ
(\mathbf K^{A_1,\mathbf{i}})^T
\mathbf{Y}
\mathbf K^{A_2,\mathbf{j}}
\mathbf{j}
\nonumber\\
&\quad
-
2(\mathbf K^{A_1,\mathbf{i}})^T
\mathbf{X}\mathbf S
\circ
\mathbf C
\mathbf{Y}
\mathbf K^{A_2,\mathbf{j}},
\\
\mathcal C_{11}(\mathbf{X},\mathbf{Y})
&=
-\mathbf{i}
(\mathbf K^{A_1,\mathbf{i}})^T
\mathbf{X}
\overline{\mathbf K}^{A_2,\mathbf{j}}
\circ
(\mathbf K^{A_1,\mathbf{i}})^T
\mathbf{Y}
\mathbf K^{A_2,\mathbf{j}}
\mathbf{j}
\nonumber\\
&\quad
+
2\mathbf{i}
(\mathbf K^{A_1,\mathbf{i}})^T
\mathbf{X}\mathbf S
\circ
\mathbf C
\mathbf{Y}
\mathbf K^{A_2,\mathbf{j}}.
\end{align*}
The operators $\mathcal C_{00}, \mathcal C_{01}, \mathcal C_{10},
\mathcal C_{11}$ differ in the placement of the quaternion units
$\mathbf{i}$ and $\mathbf{j}$ and in the signs of the correction
terms. Here $\mathbf C$ and $\mathbf S$ are defined entrywise by
\begin{align*}
[\mathbf C]_{m_1,y_1}
&=
\frac{1}{\sqrt{N_1}}
\cos\!\left(
\frac{a_1}{2b_1}y_1^2\Delta t_1^2
-\frac{2\pi}{N_1}m_1y_1
+\frac{d_1}{2b_1}m_1^2\Delta u_1^2
\right),
\\
[\mathbf S]_{z_2,m_2}
&=
\frac{1}{\sqrt{N_2}}
\sin\!\left(
\frac{a_2}{2b_2}z_2^2\Delta t_2^2
-\frac{2\pi}{N_2}m_2z_2
+\frac{d_2}{2b_2}m_2^2\Delta u_2^2
\right).
\end{align*}
Multiplication by $\mathbf{i}$ and $\mathbf{j}$ is entrywise from
the left and right, respectively.

\begin{theorem}[Spatial convolution theorem]
Let
$\widehat{\mathbf{X}},\widehat{\mathbf{Y}}\in\mathbb{DH}^{N_1\times N_2}$. Then
\begin{align}
&\mathcal L^{TS}_{A_1,\mathbf{i},A_2,\mathbf{j}}
[\widehat{\mathbf{X}}*\widehat{\mathbf{Y}}](m_1,m_2)
\nonumber\\
&=
\sqrt{N_1N_2}\,
 E^{\mathbf{i}}_{m_1}
\Big[
\mathcal C_{00}(\mathbf{X}_{i1},\mathbf{Y}_{j1})
+
\mathcal C_{01}(\mathbf{X}_{i1},\mathbf{Y}_{j2})
+
\mathcal C_{10}(\mathbf{X}_{i2},\mathbf{Y}_{j1})
+
\mathcal C_{11}(\mathbf{X}_{i2},\mathbf{Y}_{j2})
\nonumber\\
&\qquad
+\varepsilon
\Big(
\mathcal C_{00}(\mathbf{X}_{i1},\mathbf{Y}_{j3})
+
\mathcal C_{01}(\mathbf{X}_{i1},\mathbf{Y}_{j4})
+
\mathcal C_{10}(\mathbf{X}_{i2},\mathbf{Y}_{j3})
+
\mathcal C_{11}(\mathbf{X}_{i2},\mathbf{Y}_{j4})
\nonumber\\
&\qquad\qquad
+
\mathcal C_{00}(\mathbf{X}_{i3},\mathbf{Y}_{j1})
+
\mathcal C_{01}(\mathbf{X}_{i3},\mathbf{Y}_{j2})
+
\mathcal C_{10}(\mathbf{X}_{i4},\mathbf{Y}_{j1})
+
\mathcal C_{11}(\mathbf{X}_{i4},\mathbf{Y}_{j2})
\Big)
\Big]_{m_1,m_2}
E^{\mathbf{j}}_{m_2},
\label{eq:spatial_convolution_theorem}
\end{align}
where
$
E^{\mathbf{i}}_{m_1}
=
\exp\!\left(
-\mathbf{i}\frac{d_1}{2b_1}
m_1^2\Delta u_1^2
\right),
 E^{\mathbf{j}}_{m_2}
=
\exp\!\left(
-\mathbf{j}\frac{d_2}{2b_2}
m_2^2\Delta u_2^2
\right).$
\end{theorem}

\begin{proof}
Substituting \eqref{eq:spatial_convolution} into the TDQ-DLCT and
using $y_r=[n_r-z_r]_{N_r}$, $r=1,2$, the quadratic terms in $n_r$
cancel. For fixed $z_r$, the map $n_r\mapsto y_r$ is a permutation of
$\{0,\ldots,N_r-1\}$, so the sums may be rewritten independently over
$z_r$ and $y_r$.

Moreover, for some integer $q_r$,
$z_r+y_r=n_r+q_rN_r,$
and hence
\[
e^{-\mathbf{i}\frac{2\pi}{N_1}m_1n_1}
=
e^{-\mathbf{i}\frac{2\pi}{N_1}m_1(z_1+y_1)},
\qquad
e^{-\mathbf{j}\frac{2\pi}{N_2}m_2n_2}
=
e^{-\mathbf{j}\frac{2\pi}{N_2}m_2(z_2+y_2)}.
\]
Therefore,
\begin{align}
&\mathcal L^{TS}_{A_1,\mathbf{i},A_2,\mathbf{j}}
[\widehat{\mathbf{X}}*\widehat{\mathbf{Y}}](m_1,m_2)
\nonumber\\
&=
\frac{1}{\sqrt{N_1N_2}}
\sum_{z_1,z_2,y_1,y_2}
e^{\mathbf{i}\left(
\frac{a_1}{2b_1}z_1^2\Delta t_1^2
-\frac{2\pi}{N_1}m_1z_1
\right)}
\nonumber\\
&\quad
e^{\mathbf{i}\left(
\frac{a_1}{2b_1}y_1^2\Delta t_1^2
-\frac{2\pi}{N_1}m_1y_1
+\frac{d_1}{2b_1}m_1^2\Delta u_1^2
\right)}
[\widehat{\mathbf{X}}]_{z_1,z_2}
[\widehat{\mathbf{Y}}]_{y_1,y_2}
\nonumber\\
&\quad
e^{\mathbf{j}\left(
\frac{a_2}{2b_2}y_2^2\Delta t_2^2
-\frac{2\pi}{N_2}m_2y_2
+\frac{d_2}{2b_2}m_2^2\Delta u_2^2
\right)}
\nonumber\\
&\quad
e^{\mathbf{j}\left(
\frac{a_2}{2b_2}z_2^2\Delta t_2^2
-\frac{2\pi}{N_2}m_2z_2
\right)}.
\label{eq:spatial_proof_core}
\end{align}

Using the decompositions
\begin{align*}
\widehat{\mathbf{X}}
&=
\mathbf{X}_{i1}
+
\mathbf{X}_{i2}\mathbf{j}
+
\varepsilon
\left(
\mathbf{X}_{i3}
+
\mathbf{X}_{i4}\mathbf{j}
\right),
\\
\widehat{\mathbf{Y}}
&=
\mathbf{Y}_{j1}
+
\mathbf{i} \mathbf{Y}_{j2}
+
\varepsilon
\left(
\mathbf{Y}_{j3}
+
\mathbf{i} \mathbf{Y}_{j4}
\right),
\end{align*}
we consider the representative primal term
$[\mathbf{X}_{i1}]_{z_1,z_2}[\mathbf{Y}_{j1}]_{y_1,y_2}$.
For notational convenience, define
\begin{align*}
\theta_1(r;m_1)
&=
\frac{a_1}{2b_1}r^2\Delta t_1^2
-\frac{2\pi}{N_1}m_1r
+\frac{d_1}{2b_1}m_1^2\Delta u_1^2,
\\
\theta_2(r;m_2)
&=
\frac{a_2}{2b_2}r^2\Delta t_2^2
-\frac{2\pi}{N_2}m_2r
+\frac{d_2}{2b_2}m_2^2\Delta u_2^2.
\end{align*}
After introducing the output chirp factors
$E^{\mathbf{i}}_{m_1}$ and $E^{\mathbf{j}}_{m_2}$, the contribution
of this component pair in \eqref{eq:spatial_proof_core} can be
written as
\begin{align}
T(\mathbf{X}_{i1},\mathbf{Y}_{j1})
&=
\frac{1}{\sqrt{N_1N_2}}\,
 E^{\mathbf{i}}_{m_1}
\sum_{z_1,z_2,y_1,y_2}
e^{\mathbf{i}\theta_1(z_1;m_1)}
[\mathbf{X}_{i1}]_{z_1,z_2}
\nonumber\\
&\quad
e^{\mathbf{i}\theta_1(y_1;m_1)}
e^{\mathbf{j}\theta_2(z_2;m_2)}
[\mathbf{Y}_{j1}]_{y_1,y_2}
e^{\mathbf{j}\theta_2(y_2;m_2)}
 E^{\mathbf{j}}_{m_2}.
\label{eq:spatial_T00_start}
\end{align}

The two middle phase factors cannot be interchanged directly because
they are generated by different quaternion units. For real $A$ and
$B$,
\begin{equation}
e^{\mathbf{i}A}e^{\mathbf{j}B}
=
e^{-\mathbf{j}B}e^{\mathbf{i}A}
+
2\mathbf{j}\sin(B)\cos(A).
\label{eq:kernel_commutation_identity}
\end{equation}
Equation~\eqref{eq:kernel_commutation_identity} follows from
$e^{\mathbf{i}A}e^{\mathbf{j}B}
= (\cos A + \mathbf{i}\sin A)(\cos B + \mathbf{j}\sin B)$
and the quaternion identities
$\mathbf{i}\mathbf{j} = \mathbf{k} = -\mathbf{j}\mathbf{i}$.
Applying \eqref{eq:kernel_commutation_identity} with
$A=\theta_1(y_1;m_1)$, $B=\theta_2(z_2;m_2),$
splits \eqref{eq:spatial_T00_start} into two terms,
\[
T(\mathbf{X}_{i1},\mathbf{Y}_{j1})
=
T_{00}^{(1)}+T_{00}^{(2)}.
\]
For the first term,
\begin{align*}
T_{00}^{(1)}
&=
\frac{1}{\sqrt{N_1N_2}}\,
 E^{\mathbf{i}}_{m_1}
\sum_{z_1,z_2,y_1,y_2}
e^{\mathbf{i}\theta_1(z_1;m_1)}
[\mathbf{X}_{i1}]_{z_1,z_2}
e^{-\mathbf{j}\theta_2(z_2;m_2)}
\nonumber\\
&\quad
e^{\mathbf{i}\theta_1(y_1;m_1)}
[\mathbf{Y}_{j1}]_{y_1,y_2}
e^{\mathbf{j}\theta_2(y_2;m_2)}
 E^{\mathbf{j}}_{m_2}.
\end{align*}
The independent sums over $(z_1,z_2)$ and $(y_1,y_2)$ separate,
giving
\begin{align}
T_{00}^{(1)}
=
\sqrt{N_1N_2}\,
 E^{\mathbf{i}}_{m_1}
\Big[
(\mathbf K^{A_1,\mathbf{i}})^T
\mathbf{X}_{i1}
\overline{\mathbf K}^{A_2,\mathbf{j}}
\circ
(\mathbf K^{A_1,\mathbf{i}})^T
\mathbf{Y}_{j1}
\mathbf K^{A_2,\mathbf{j}}
\Big]_{m_1,m_2}
 E^{\mathbf{j}}_{m_2}.
\label{eq:spatial_T00_first}
\end{align}
The factor $\sqrt{N_1N_2}$ appears because each of the two
TDQ-DLCT-type matrix products contains a factor
$(N_1N_2)^{-1/2}$, whereas \eqref{eq:spatial_T00_start} contains only
one such normalization.

The second term generated by
\eqref{eq:kernel_commutation_identity} is
\begin{align*}
T_{00}^{(2)}
&=
\frac{2}{\sqrt{N_1N_2}}\,
 E^{\mathbf{i}}_{m_1}
\sum_{z_1,z_2,y_1,y_2}
e^{\mathbf{i}\theta_1(z_1;m_1)}
[\mathbf{X}_{i1}]_{z_1,z_2}
\nonumber\\
&\quad
\mathbf{j}
\sin\!\big(\theta_2(z_2;m_2)\big)
\cos\!\big(\theta_1(y_1;m_1)\big)
[\mathbf{Y}_{j1}]_{y_1,y_2}
e^{\mathbf{j}\theta_2(y_2;m_2)}
 E^{\mathbf{j}}_{m_2}.
\end{align*}
Since $\mathbf{Y}_{j1}\in\mathbb C_{\mathbf{j}}$ and the $\mathbf{j}$-generated
kernel belongs to the same complex subalgebra, $\mathbf{j}$ commutes
with both of them. Thus the factor $\mathbf{j}$ may be moved to the
right. Using
\begin{align*}
[\mathbf C]_{m_1,y_1}
&=
\frac{1}{\sqrt{N_1}}
\cos\!\big(\theta_1(y_1;m_1)\big),
\\
[\mathbf S]_{z_2,m_2}
&=
\frac{1}{\sqrt{N_2}}
\sin\!\big(\theta_2(z_2;m_2)\big),
\end{align*}
we obtain
\begin{align}
T_{00}^{(2)}
=
2\sqrt{N_1N_2}\,
 E^{\mathbf{i}}_{m_1}
\Big[
(\mathbf K^{A_1,\mathbf{i}})^T
\mathbf{X}_{i1}\mathbf S
\circ
\mathbf C
\mathbf{Y}_{j1}
\mathbf K^{A_2,\mathbf{j}}
\mathbf{j}
\Big]_{m_1,m_2}
 E^{\mathbf{j}}_{m_2}.
\label{eq:spatial_T00_second}
\end{align}

Combining \eqref{eq:spatial_T00_first} and
\eqref{eq:spatial_T00_second} yields
\begin{align*}
T(\mathbf{X}_{i1},\mathbf{Y}_{j1})
&=
\sqrt{N_1N_2}\,
 E^{\mathbf{i}}_{m_1}
\Big[
(\mathbf K^{A_1,\mathbf{i}})^T
\mathbf{X}_{i1}
\overline{\mathbf K}^{A_2,\mathbf{j}}
\circ
(\mathbf K^{A_1,\mathbf{i}})^T
\mathbf{Y}_{j1}
\mathbf K^{A_2,\mathbf{j}}
\nonumber\\
&\qquad
+
2(\mathbf K^{A_1,\mathbf{i}})^T
\mathbf{X}_{i1}\mathbf S
\circ
\mathbf C
\mathbf{Y}_{j1}
\mathbf K^{A_2,\mathbf{j}}
\mathbf{j}
\Big]_{m_1,m_2}
 E^{\mathbf{j}}_{m_2}
\nonumber\\
&=
\sqrt{N_1N_2}\,
 E^{\mathbf{i}}_{m_1}
\left[
\mathcal C_{00}
(\mathbf{X}_{i1},\mathbf{Y}_{j1})
\right]_{m_1,m_2}
 E^{\mathbf{j}}_{m_2}.
\end{align*}

The remaining primal component products are obtained analogously:
\begin{align*}
T(\mathbf{X}_{i1},\mathbf{i} \mathbf{Y}_{j2})
&=
\sqrt{N_1N_2}\,
 E^{\mathbf{i}}_{m_1}
\left[
\mathcal C_{01}
(\mathbf{X}_{i1},\mathbf{Y}_{j2})
\right]_{m_1,m_2}
 E^{\mathbf{j}}_{m_2},
\\
T(\mathbf{X}_{i2}\mathbf{j},\mathbf{Y}_{j1})
&=
\sqrt{N_1N_2}\,
 E^{\mathbf{i}}_{m_1}
\left[
\mathcal C_{10}
(\mathbf{X}_{i2},\mathbf{Y}_{j1})
\right]_{m_1,m_2}
 E^{\mathbf{j}}_{m_2},
\\
T(\mathbf{X}_{i2}\mathbf{j},\mathbf{i} \mathbf{Y}_{j2})
&=
\sqrt{N_1N_2}\,
 E^{\mathbf{i}}_{m_1}
\left[
\mathcal C_{11}
(\mathbf{X}_{i2},\mathbf{Y}_{j2})
\right]_{m_1,m_2}
 E^{\mathbf{j}}_{m_2}.
\end{align*}

Finally, since $\varepsilon$ commutes with quaternion-valued
quantities and satisfies $\varepsilon^2=0$,
\[
(a+\varepsilon b)(c+\varepsilon d)
=
ac+\varepsilon(ad+bc).
\]
Thus, the dual part consists of the eight primal--dual cross terms.
Applying the four identities above and collecting all contributions
yields \eqref{eq:spatial_convolution_theorem}. The four primal terms
yield the first line of \eqref{eq:spatial_convolution_theorem}, and
the eight primal--dual cross terms yield the $\varepsilon$-part.
\end{proof}

\subsection{Frequency Convolution}
\label{subsec:frequency_convolution}

The spatial convolution of the previous subsection corresponds, in
the transform domain, to an ordered Hadamard product. Dually, we may
define a convolution directly in the transform domain and ask what it
corresponds to in the spatial domain. The following definition
formalizes this dual construction.

\begin{definition}[Frequency convolution]
\label{fre_con}
Let $\widehat{\mathbf{X}},\widehat{\mathbf{Y}}\in\mathbb{DH}^{N_1\times N_2}$.
Their TDQ-DLCT-induced frequency convolution is defined by
\begin{equation}
\widehat{\mathbf{X}}
\circledast_{A_1,A_2}
\widehat{\mathbf{Y}}
=
\left(
\mathcal L_{A_1,\mathbf{i},A_2,\mathbf{j}}^{TS}
\right)^{-1}
\left[
\mathcal L_{A_1,\mathbf{i},A_2,\mathbf{j}}^{TS}
[\widehat{\mathbf{X}}]
\circ
\mathcal L_{A_1,\mathbf{i},A_2,\mathbf{j}}^{TS}
[\widehat{\mathbf{Y}}]
\right].
\label{eq:frequency_convolution}
\end{equation}
\end{definition}

By construction,
\begin{equation*}
\mathcal L_{A_1,\mathbf{i},A_2,\mathbf{j}}^{TS}
\left[
\widehat{\mathbf{X}}
\circledast_{A_1,A_2}
\widehat{\mathbf{Y}}
\right]
=
\mathcal L_{A_1,\mathbf{i},A_2,\mathbf{j}}^{TS}
[\widehat{\mathbf{X}}]
\circ
\mathcal L_{A_1,\mathbf{i},A_2,\mathbf{j}}^{TS}
[\widehat{\mathbf{Y}}].
\end{equation*}

\begin{remark}[Duality of the two convolutions]
\label{rem:convolution_duality}
The spatial convolution admits a structured
transform-domain representation involving ordered
Hadamard products and noncommutative correction terms,
as shown in \eqref{eq:spatial_convolution_theorem}.
In contrast, the frequency convolution is defined
such that its TDQ-DLCT equals the ordered Hadamard
product of the individual transforms.
\end{remark}

For
$\mathbf{X}\in\mathbb C_{\mathbf{i}}^{N_1\times N_2}$ and
$\mathbf{Y}\in\mathbb C_{\mathbf{j}}^{N_1\times N_2}$, define
\begin{align*}
\mathcal D_{00}^{A_1,A_2}(\mathbf{X},\mathbf{Y})
&=
(\mathbf K^{A_1,\mathbf{i}})^T
\mathbf{X}
\mathbf K^{A_2,\mathbf{j}}
\circ
(\mathbf K^{A_1,\mathbf{i}})^T
\mathbf{Y}
\mathbf K^{A_2,\mathbf{j}},
\\
\mathcal D_{01}^{A_1,A_2}(\mathbf{X},\mathbf{Y})
&=
\mathbf{i}
(\mathbf K^{A_1,\mathbf{i}})^T
\mathbf{X}
\overline{\mathbf K}^{A_2,\mathbf{j}}
\circ
(\mathbf K^{A_1,\mathbf{i}})^T
\mathbf{Y}
\mathbf K^{A_2,\mathbf{j}},
\\
\mathcal D_{10}^{A_1,A_2}(\mathbf{X},\mathbf{Y})
&=
(\mathbf K^{A_1,\mathbf{i}})^T
\mathbf{X}
\mathbf K^{A_2,\mathbf{j}}
\circ
(\overline{\mathbf K}^{A_1,\mathbf{i}})^T
\mathbf{Y}
\mathbf K^{A_2,\mathbf{j}}
\mathbf{j},
\\
\mathcal D_{11}^{A_1,A_2}(\mathbf{X},\mathbf{Y})
&=
-\mathbf{i}
(\mathbf K^{A_1,\mathbf{i}})^T
\mathbf{X}
\overline{\mathbf K}^{A_2,\mathbf{j}}
\circ
(\overline{\mathbf K}^{A_1,\mathbf{i}})^T
\mathbf{Y}
\mathbf K^{A_2,\mathbf{j}}
\mathbf{j}.
\end{align*}
These four operators mirror the spatial operators
$\mathcal C_{00}, \mathcal C_{01}, \mathcal C_{10}, \mathcal C_{11}$
in the transform domain and account for the same four
quaternion-component pairings; they differ in that the canonical
kernels appear without their input-chirp and output-chirp
renormalization, since the convolution is already defined in the
transform domain.

\begin{theorem}[Frequency convolution theorem]
Let $\widehat{\mathbf{X}},\widehat{\mathbf{Y}}\in\mathbb{DH}^{N_1\times N_2}$. Then
\begin{align}
&\mathcal L_{A_1,\mathbf{i},A_2,\mathbf{j}}^{TS}
\left[
\widehat{\mathbf{X}}
\circledast_{A_1,A_2}
\widehat{\mathbf{Y}}
\right]
\nonumber\\
&=
\mathcal D_{00}^{A_1,A_2}
(\mathbf{X}_{i1},\mathbf{Y}_{j1})
+
\mathcal D_{01}^{A_1,A_2}
(\mathbf{X}_{i1},\mathbf{Y}_{j2})
+
\mathcal D_{10}^{A_1,A_2}
(\mathbf{X}_{i2},\mathbf{Y}_{j1})
+
\mathcal D_{11}^{A_1,A_2}
(\mathbf{X}_{i2},\mathbf{Y}_{j2})
\nonumber\\
&\quad
+\varepsilon
\Big[
\mathcal D_{00}^{A_1,A_2}
(\mathbf{X}_{i1},\mathbf{Y}_{j3})
+
\mathcal D_{01}^{A_1,A_2}
(\mathbf{X}_{i1},\mathbf{Y}_{j4})
+
\mathcal D_{10}^{A_1,A_2}
(\mathbf{X}_{i2},\mathbf{Y}_{j3})
+
\mathcal D_{11}^{A_1,A_2}
(\mathbf{X}_{i2},\mathbf{Y}_{j4})
\nonumber\\
&\qquad
+
\mathcal D_{00}^{A_1,A_2}
(\mathbf{X}_{i3},\mathbf{Y}_{j1})
+
\mathcal D_{01}^{A_1,A_2}
(\mathbf{X}_{i3},\mathbf{Y}_{j2})
+
\mathcal D_{10}^{A_1,A_2}
(\mathbf{X}_{i4},\mathbf{Y}_{j1})
+
\mathcal D_{11}^{A_1,A_2}
(\mathbf{X}_{i4},\mathbf{Y}_{j2})
\Big].
\label{eq:frequency_convolution_theorem}
\end{align}
\end{theorem}

\begin{proof}
By Definition~\ref{fre_con},
\begin{equation}
\mathcal L_{A_1,\mathbf{i},A_2,\mathbf{j}}^{TS}
\left[
\widehat{\mathbf{X}}
\circledast_{A_1,A_2}
\widehat{\mathbf{Y}}
\right]
=
\mathcal L_{A_1,\mathbf{i},A_2,\mathbf{j}}^{TS}
[\widehat{\mathbf{X}}]
\circ
\mathcal L_{A_1,\mathbf{i},A_2,\mathbf{j}}^{TS}
[\widehat{\mathbf{Y}}].
\label{eq:frequency_proof_start}
\end{equation}

Using the decompositions
$
\widehat{\mathbf{X}}
=
\mathbf{X}_{i1}
+
\mathbf{X}_{i2}\mathbf{j}
+
\varepsilon
\left(
\mathbf{X}_{i3}
+
\mathbf{X}_{i4}\mathbf{j}
\right),$
$\widehat{\mathbf{Y}}
=
\mathbf{Y}_{j1}
+
\mathbf{i} \mathbf{Y}_{j2}
+
\varepsilon
\left(
\mathbf{Y}_{j3}
+
\mathbf{i} \mathbf{Y}_{j4}
\right),$
the ordered Hadamard product on the right-hand side of
\eqref{eq:frequency_proof_start} can be expanded componentwise.

We consider the interaction between
$\mathbf{X}_{i1}$ and $\mathbf{i} \mathbf{Y}_{j2}$ as a representative term. Its
contribution is
\begin{equation}
\left[
(\mathbf K^{A_1,\mathbf{i}})^T
\mathbf{X}_{i1}
\mathbf K^{A_2,\mathbf{j}}
\right]_{m_1,m_2}
\mathbf{i}
\left[
(\mathbf K^{A_1,\mathbf{i}})^T
\mathbf{Y}_{j2}
\mathbf K^{A_2,\mathbf{j}}
\right]_{m_1,m_2}.
\label{eq:frequency_D01_start}
\end{equation}

Using
$
[\mathbf K^{A_2,\mathbf{j}}]_{n_2,m_2}\mathbf{i}
=
\mathbf{i}
[\overline{\mathbf K}^{A_2,\mathbf{j}}]_{n_2,m_2},$
which follows from
$\mathbf{j}\mathbf{i} = -\mathbf{i}\mathbf{j}$ and the entrywise
form of $\mathbf K^{A_2,\mathbf{j}}$, and the fact that
$[\mathbf K^{A_1,\mathbf{i}}]_{n_1,m_1}$ and
$\mathbf{X}_{i1}(n_1,n_2)$ belong to $\mathbb C_{\mathbf{i}}$, we obtain
\begin{equation}
\left[
(\mathbf K^{A_1,\mathbf{i}})^T
\mathbf{X}_{i1}
\mathbf K^{A_2,\mathbf{j}}
\right]_{m_1,m_2}
\mathbf{i}
=
\mathbf{i}
\left[
(\mathbf K^{A_1,\mathbf{i}})^T
\mathbf{X}_{i1}
\overline{\mathbf K}^{A_2,\mathbf{j}}
\right]_{m_1,m_2}.
\label{eq:frequency_D01_commutation}
\end{equation}

Substituting \eqref{eq:frequency_D01_commutation} into
\eqref{eq:frequency_D01_start} gives
\begin{align*}
&
\mathbf{i}
\left[
(\mathbf K^{A_1,\mathbf{i}})^T
\mathbf{X}_{i1}
\overline{\mathbf K}^{A_2,\mathbf{j}}
\right]_{m_1,m_2}
\left[
(\mathbf K^{A_1,\mathbf{i}})^T
\mathbf{Y}_{j2}
\mathbf K^{A_2,\mathbf{j}}
\right]_{m_1,m_2}
\nonumber\\
&=
\left[
\mathbf{i}
(\mathbf K^{A_1,\mathbf{i}})^T
\mathbf{X}_{i1}
\overline{\mathbf K}^{A_2,\mathbf{j}}
\circ
(\mathbf K^{A_1,\mathbf{i}})^T
\mathbf{Y}_{j2}
\mathbf K^{A_2,\mathbf{j}}
\right]_{m_1,m_2}
\nonumber\\
&=
\left[
\mathcal D_{01}^{A_1,A_2}
(\mathbf{X}_{i1},\mathbf{Y}_{j2})
\right]_{m_1,m_2}.
\end{align*}
Thus, the component pair $\mathbf{X}_{i1}$ and $\mathbf{i} \mathbf{Y}_{j2}$ produces
precisely $\mathcal D_{01}^{A_1,A_2}(\mathbf{X}_{i1},\mathbf{Y}_{j2}).$

The remaining component pairs follow analogously. In particular,
\[
\mathbf{j}
[\mathbf K^{A_1,\mathbf{i}}]_{n_1,m_1}
=
[\overline{\mathbf K}^{A_1,\mathbf{i}}]_{n_1,m_1}
\mathbf{j},
\qquad
\mathbf{j}\mathbf{i}=-\mathbf{i}\mathbf{j},
\]
lead to the operators
$\mathcal D_{10}^{A_1,A_2}$ and
$\mathcal D_{11}^{A_1,A_2}$, while the component pair without
basis-unit interaction gives $\mathcal D_{00}^{A_1,A_2}$.

Finally, since $\varepsilon$ commutes with quaternion-valued
quantities and satisfies $\varepsilon^2=0$,
$(a+\varepsilon b)(c+\varepsilon d)
=ac+\varepsilon(ad+bc).$
Hence, the dual part consists of the eight primal--dual cross terms.
Collecting all primal and dual contributions yields
\eqref{eq:frequency_convolution_theorem}, which completes the proof.
\end{proof}

\FloatBarrier

\section{Experiments}\label{sec:experiments}

To assess the practical potential of the proposed theoretical
framework, we first evaluate the fast DQ-DLCT, which supports the
efficient implementation of the subsequent denoising method. We then
introduce a simple chirp-adaptive DQ-DLCT algorithm and evaluate its
denoising performance.

\subsection{Computational Efficiency}\label{subsec:fast}

We compare direct evaluation of \eqref{eq:2d_two_sided_DQ-DLCT}
with the fast algorithm based on
\eqref{eq:2d_DQ-DLCT_factorization}. The direct method requires
$\mathcal O(N^4)$ operations on an $N\times N$ array, whereas the
fast method requires $\mathcal O(N^2\log N)$ operations.

For each size $N$, we generate an $N\times N$ dual-quaternion array
whose eight real coefficients at each location are independently
drawn from the standard normal distribution. The same array is
supplied to both implementations. We use random seed zero,
$\Delta t_1=\Delta t_2=0.1$, and
\[
A_1=A_2=
\begin{bmatrix}
1&1\\
0&1
\end{bmatrix},
\qquad
\Delta u_s=\frac{2\pi b_s}{N\Delta t_s},
\quad s\in\{1,2\}.
\]
This matrix satisfies $a_sd_s-b_sc_s=1$ and $b_s\ne0$.

We test powers of two from $N=1$ to $N=512$, with direct evaluation
measured up to $N=256$. Input generation is excluded from timing.
The fast method is timed over ten calls at each size. The direct
method is timed over ten calls for $N\le32$, three calls for $N=64$,
and one call for $N=128$ and $N=256$. Reported runtimes are
arithmetic means of the individual wall-clock durations.

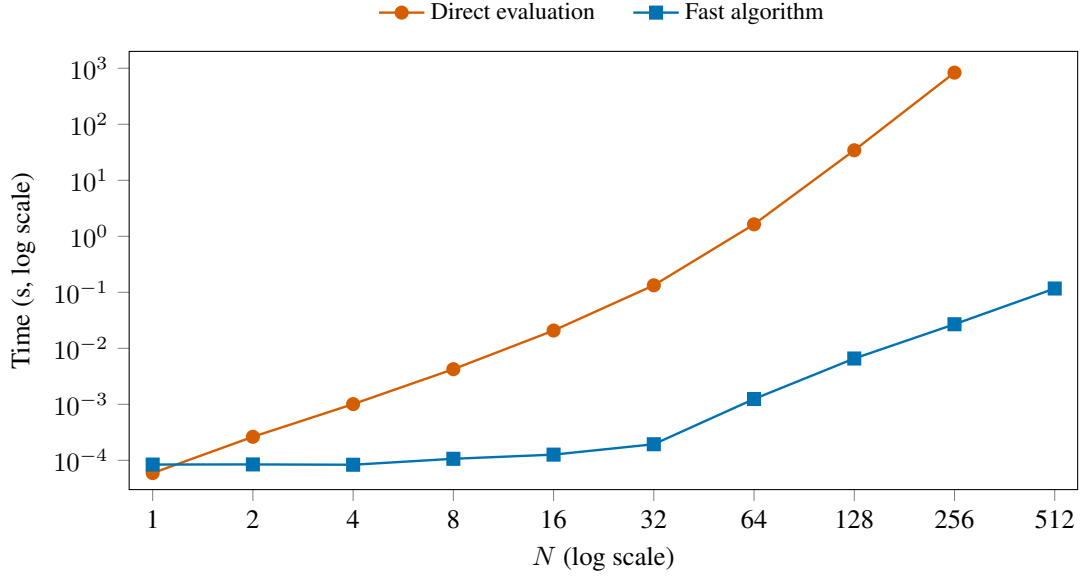
\begin{figure}[t]
\centering
\begin{tikzpicture}
\begin{axis}[
  width=0.76\linewidth, height=58mm, scale only axis,
  xlabel={$N$ (log scale)}, ylabel={Time (s, log scale)},
  xmode=log, log basis x=2, xmin=0.85, xmax=600,
  xtick={1,2,4,8,16,32,64,128,256,512},
  xticklabels={1,2,4,8,16,32,64,128,256,512},
  ymode=log, ymin=0.00003, ymax=2000,
  ytick={0.0001,0.001,0.01,0.1,1,10,100,1000},
  grid=none, minor tick length=0pt,
  minor tick style={draw=none},
  tick align=outside, tick pos=left,
  legend style={draw=none, fill=none, font=\small,
                at={(0.5,1.04)}, anchor=south,
                legend columns=2,
                /tikz/every even column/.append style={column sep=12pt}},
  every axis plot/.append style={line width=0.9pt, mark size=2.2pt},
]
\addplot[color={rgb,255:red,213;green,94;blue,0}, mark=*]
  coordinates {(1,0.0000592199969105) (2,0.000262949999887)
               (4,0.00100784999086) (8,0.00423998000915)
               (16,0.0207564299926) (32,0.133262980013)
               (64,1.63456056667) (128,34.3115529)
               (256,835.0887012)};
\addlegendentry{Direct evaluation}
\addplot[color={rgb,255:red,0;green,114;blue,178}, mark=square*]
  coordinates {(1,0.000084020011127) (2,0.0000844600028358)
               (4,0.0000833799887914) (8,0.000106469995808)
               (16,0.000126379996073) (32,0.000193730002502)
               (64,0.0012431900017) (128,0.00657164001605)
               (256,0.0268850200053) (512,0.11684525999)};
\addlegendentry{Fast algorithm}
\end{axis}
\end{tikzpicture}
\caption{Measured runtimes of direct evaluation and the fast
two-dimensional TDQ-DLCT for the specified canonical matrix.}
\label{fig:runtime}
\end{figure}

\begin{figure}[t]
\centering
\begin{tikzpicture}
\begin{axis}[
  width=0.7\linewidth, height=55mm,
  xlabel={$N$ (log scale)},
  ylabel={$E_{\infty}(N)$ (log scale)},
  xmode=log, log basis x=2, xmin=1.7, xmax=300,
  xtick={2,4,8,16,32,64,128,256},
  ymode=log, ymin=5e-15, ymax=1e-12,
  grid=both, grid style={gray!20},
  every axis plot/.append style={line width=0.9pt, mark size=2.2pt},
]
\addplot[color={rgb,255:red,213;green,94;blue,0}, mark=*]
  coordinates {(2,7.327471962526033e-15)
               (4,2.0117241206207837e-13)
               (8,1.9501067427540875e-13)
               (16,2.892130979148533e-13)
               (32,3.2429614549300823e-13)
               (64,3.865796571744795e-13)
               (128,3.788080960021034e-13)
               (256,5.129230373768223e-13)};
\end{axis}
\end{tikzpicture}
\caption{Measured maximum coefficientwise discrepancy between the
fast and direct TDQ-DLCT. The zero discrepancy at $N=1$ cannot be
shown on logarithmic axes.}
\label{fig:error_scaling}
\end{figure}
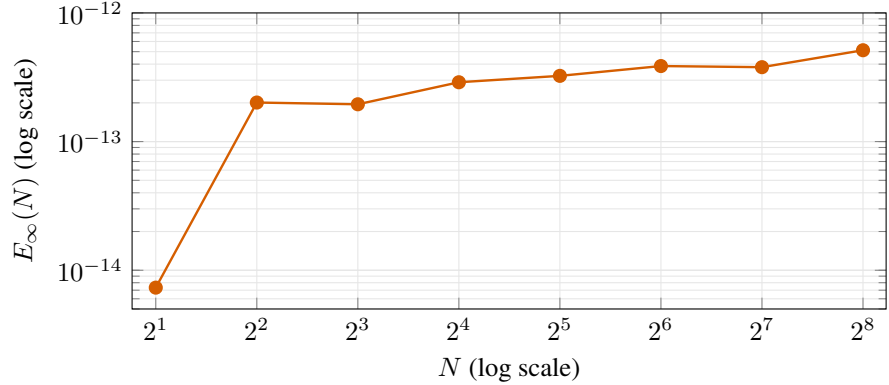

For every size at which both implementations are evaluated, we
compare all eight real coefficients of their complete outputs.
The maximum coefficientwise discrepancy is
\[
E_{\infty}(N)
=
\max_{\substack{0\le m_1,m_2<N\\1\le r\le8}}
\left|
\widehat F_{\mathrm{fast},r}(m_1,m_2)
-
\widehat F_{\mathrm{direct},r}(m_1,m_2)
\right|.
\]
Figure~\ref{fig:error_scaling} displays the measured discrepancies.
The value at $N=1$ is exactly zero and is omitted from the
logarithmic plot. Across all tested sizes, the largest discrepancy
is $5.13\times10^{-13}$, attained at $N=256$.

As shown in Figure~\ref{fig:runtime}, the fast method is slightly
slower at $N=1$ but faster from $N=2$ onward. At $N=256$, direct
evaluation takes approximately $835.09$~s, compared with
$0.0269$~s for the fast algorithm, giving a measured speedup of
about $3.11\times10^4$. Least-squares fits of log runtime against
log $N$ over $N=16,\ldots,256$ for the direct method and
$N=16,\ldots,512$ for the fast method give slopes of approximately
$3.86$ and $2.09$, respectively. These finite-range measurements
are consistent with the respective complexity estimates, although
the fitted slopes and speedup depend on the implementation and
computing environment. The small discrepancies in
Figure~\ref{fig:error_scaling} confirm numerical agreement between
the direct and fast implementations for the tested transform.
\FloatBarrier

\subsection{Denoising}\label{subsec:denoise}

\subsubsection{Algorithm}\label{subsec:filter}

We introduce chirp-adaptive DQ-DLCT denoising (CA-DQ-DLCT) with the
same two-stage procedure for pose sequences and fields. Spectral
filtering first estimates rotation, then refines position using that
estimate. Chirp rates adapt to observations; window and gain settings
are fitted before testing, and the output is a unit dual quaternion.

Output chirps preserve coefficient powers and cancel in the inverse
transform when filtering uses a shared real gain based on those
powers. We therefore set $d_s=0$ and absorb the input phase
coefficient $a_s\Delta t_s^2/(2b_s)$ into an effective chirp rate
$\alpha_s$. Using normalized sampling coordinates, we obtain the
equivalent parameterization
\begin{equation*}
A_s=A(\alpha_s)=\begin{bmatrix} \alpha_s & 1\\ -1 & 0\end{bmatrix},
\qquad \Delta t_s=\Delta u_s=\sqrt{2\pi/N_s},
\qquad s\in\{1,2\},
\end{equation*}
where $N_s$ is the number of samples along axis $s$. Setting the
active chirp rates to zero recovers the corresponding dual quaternion
Fourier transform.

For rotation about the kernel axis with constant angular acceleration
$\gamma_s$, chirp matching \cite{xia2000chirp} gives the
phase-cancelling rate under our normalization,
\[
\alpha_s^{\mathrm m}=-\frac{\gamma_s T_s^2N_s}{4\pi},
\]
where $T_s$ is the physical sampling period. This relation motivates
chirp adaptation for accelerating motion; the algorithm selects rates
from the observations without assuming known acceleration.

We use separable square-root Hann windows and weighted overlap-add
\cite{smith2011spectral}, with half-window overlap and reflected
padding where windows are shorter than the field. In a window of size
$L_s$, the local rate is $\alpha_s L_s/N_s$. For its transform
$F_{\boldsymbol\alpha}$, define
\begin{equation}
\begin{aligned}
P_{\boldsymbol\alpha}(k)&=\sum_{b=0}^1\omega_b
 \|F_{\boldsymbol\alpha,b}(k)\|^2,
\qquad p_{\boldsymbol\alpha}(k)=\frac{P_{\boldsymbol\alpha}(k)}
 {\sum_\ell P_{\boldsymbol\alpha}(\ell)},\\
H_3(p_{\boldsymbol\alpha})&=-\tfrac12\log\sum_k p_{\boldsymbol\alpha}(k)^3,
\end{aligned}
\label{eq:entropy}
\end{equation}
where $b$ indexes the primal and dual parts and $k$ is a frequency
index. Adaptive candidates reduce this R\'enyi entropy
\cite{baraniuk2001measuring} by alternating coarse-to-fine searches
over the active axes \cite{wright2015coordinate}; zero-rate candidates
are retained.

Inverse-variance weights $\omega_b\propto v_b^{-1}$ balance the two
parts in the initial pose estimates. We estimate $v_b$ by summing
squared, centred second-difference MAD scales over its four
coefficients, divided by $6(0.67449)^2$
\cite{donoho1994ideal}. We average these estimates over the active
axes. This is an input-based noise approximation.

Following spectral subtraction \cite{boll1979} and Wiener-form
shrinkage \cite{lim1979enhancement}, each branch uses the shared real
gain
\begin{equation}
G(k)=\frac{(P(k)-c\nu)_+}{P(k)},\qquad
\nu=\Big(\sum_b\omega_bv_b\Big)\frac{\|w\|^2}{\prod_s L_s},
\label{eq:gain}
\end{equation}
with $G=0$ for a zero denominator. Both denoising experiments use the
unsmoothed power $P$ and reconstruct by weighted overlap-add.

Both experiments use $c\in\{0.25,0.5,1,2,4,8,16,32\}$, zero or
entropy-selected rates, and an input fallback. Initial windows have
lengths $64$, $128$, or $512$ for sequences and sizes $64\times128$,
$128\times256$, or $256\times512$ for fields.

The rotation settings minimize fitting rotation MSE. With this
estimate fixed, position settings minimize fitting position MSE.
Validation refits both stages, excluding the held recording and every
field sharing either source; all settings are frozen before test
evaluation.

In both dimensions, we re-encode observed translation $\tilde t$
using the denoised rotation $q_R$ before refining position:
\begin{equation*}
g_P=q_R+\frac{\varepsilon}{2}(0,\tilde t)q_R,
\qquad \tilde t=2\operatorname{Vec}(\tilde d\tilde q^*).
\end{equation*}
This removes direct modulation by the observed rotation noise.
A DQ-DLCT over the full input filters both parts with weights $(0,1)$
in \eqref{eq:entropy}--\eqref{eq:gain}. We decode position after
projection, retain $q_R$, and keep the initial pose estimate as a
fallback.

Projection onto the unit dual quaternion manifold is defined by
\begin{equation*}
\Pi(\widehat f)
=
\frac{q}{\|q\|}
+\varepsilon\!\left(
\frac{d}{\|q\|}
-\frac{q}{\|q\|}
\Big\langle
\frac{q}{\|q\|},\frac{d}{\|q\|}
\Big\rangle
\right),
\qquad
\widehat f = q+\varepsilon d,
\end{equation*}
which enforces $\|q\|=1$ and $\langle q,d\rangle=0$. Uncovered samples
or vanishing primal norms retain the branch input. The final pose
combines the two branch estimates as
\begin{equation*}
\hat f=q_R+\frac{\varepsilon}{2}(0,t_P)q_R,
\qquad t_P=2\operatorname{Vec}(d_Pq_P^*).
\end{equation*}
Thus, changing the rotation estimate preserves the position branch's
decoded translation in both sequence and field denoising.

\begin{algorithm}[t]
\caption{Chirp-adaptive DQ-DLCT denoising (CA-DQ-DLCT)}
\label{alg:denoise}
\renewcommand{\baselinestretch}{1}
\small
\begin{algorithmic}[1]
\Require $\tilde f=\tilde q+\varepsilon\tilde d$; rate grids; fitted $\theta_R,\theta_P$
  ($\boldsymbol L_P=\boldsymbol N$)
\Ensure denoised unit pose field $\hat f$
\Statex $\omega_{R,j}\propto(v_{0,j}^{-1},v_{1,j}^{-1})$,
  $\omega_{P,j}=(0,1)$; $G_j$ from \eqref{eq:gain}.
\Statex $\operatorname{Rate}_{\theta}$: zero rate or entropy search
  \eqref{eq:entropy}; $\Pi_g$: unit-pose projection.
\State $\tilde t\gets2\operatorname{Vec}(\tilde d\tilde q^*)$
\For{$B\in(R,P)$ in this order}
  \State $(g_B,h_B)\gets
    \begin{cases}
      (\tilde f,\tilde f), & B=R,\\
      (q_R+\tfrac{\varepsilon}{2}(0,\tilde t)q_R,\,f_R), & B=P
    \end{cases}$
  \State \textbf{if} $\theta_B=\mathrm{id}$ \textbf{then}
    $f_B\gets h_B=q_B+\varepsilon d_B$; \textbf{continue}
  \For{each window $g_{B,j}$ specified by $\theta_B$}
    \State $\boldsymbol\alpha_j\gets
      \operatorname{Rate}_{\theta_B}(w_jg_{B,j},\omega_{B,j})$
    \State $F_j\gets\mathcal L_{\boldsymbol\alpha_j}(w_jg_{B,j})$
    \State $y_j\gets\mathcal L_{\boldsymbol\alpha_j}^{-1}(G_jF_j)$
  \EndFor
  \State $f_B\gets\Pi_{g_B}\operatorname{WOLA}_{\theta_B}\{y_j\}
    =q_B+\varepsilon d_B$
\EndFor
\State $t_P\gets2\operatorname{Vec}(d_Pq_P^*)$
\State \Return $\hat f=q_R+\tfrac{\varepsilon}{2}(0,t_P)q_R$
\end{algorithmic}
\end{algorithm}
\FloatBarrier

\subsubsection{One-dimensional Denoising}\label{subsec:denoise1d}

We evaluate the denoising method on measured pose sequences from the
Mini Wheelbot dataset \cite{hose2026wheelbot}. We extract one
$512$-sample window from each of $49$ recordings and use $24$
recordings for fitting and $25$ for testing.\footnote{Local position
values are interpreted as millimetres, consistent with native Vicon
units \cite{vicon2026sdk} and encoder displacement scales, despite the
dataset's metre label.} Before relative-pose encoding, the measured
unit quaternions are made sign-continuous. Specifically, for
consecutive samples we apply
$q_k\leftarrow -q_k$ if $\langle q_{k-1},q_k \rangle<0.$
The dual part is constructed after this alignment, so the primal and
dual components use a consistent representative of the same rigid
pose. Orientation $q$ and position $t$, relative to the first pose,
are encoded as $q+\varepsilon tq/2$, with nonzero dual parts.

All methods receive identical added Gaussian rotation-vector noise
with per-coordinate $\sigma\in\{3^\circ,10^\circ,20^\circ\}$ and
position noise of $5$~mm. We set $\alpha_1=0$ and use
$\mathcal A_0=\{-0.1+0.0025r:r=0,\ldots,80\}$ for $\alpha_2$.
Branch settings are fitted before test evaluation, with measured
trajectories as the references.

DFT \cite{smith2011spectral} and DLCT \cite{healy2016linear} filter
the eight real coefficients independently; DFT
\cite{hitzer2021introduction} and QLCT \cite{bahri2021discrete}
filter quaternion rotations and Cartesian positions separately. We
combine these transforms with a full-length taper, MAD-based spectral
subtraction \cite{boll1979}, and prescribed LCT rate $\alpha_2=0.0025$.
We also apply cubic, $11$-sample Savitzky--Golay smoothing
\cite{savitzky1964} and constant-acceleration RTS smoothing
\cite{rauch1965}, with jerk-to-measurement variance ratio selected
from $\{10^{-4},10^{-3},10^{-2},10^{-1}\}$ using noise-normalized
fitting errors. One ratio is shared across noise levels and axes.
All outputs are projected onto valid unit poses.

CA-DQ-DLCT yields the lowest pooled errors in Table~\ref{tab:denoise}.
Five further noise draws preserve its advantage over DFT/QLCT on the
same recordings. These comparisons concern the specified pipelines,
not an intrinsic advantage of the dual-quaternion representation.

\begin{table}[t]
\centering
\caption{Pooled rotation and position RMSE on $25$ held-out Mini
Wheelbot recordings, evaluated on the central $256$ samples per
window. Position noise is $5$~mm per coordinate. DLCT and QLCT use
$\alpha_2=0.0025$. Bold identifies the lowest error in each column.}
\label{tab:denoise}
\vspace{6pt}
\sisetup{table-format=2.3}
\small
\setlength{\tabcolsep}{3pt}
\begin{tabular}{l SSS SSS}
\toprule
& \multicolumn{3}{c}{Rotation RMSE (deg)} & \multicolumn{3}{c}{Position RMSE (mm)} \\
\cmidrule(lr){2-4}\cmidrule(lr){5-7}
Method & {$3^\circ$} & {$10^\circ$} & {$20^\circ$} & {$3^\circ$} & {$10^\circ$} & {$20^\circ$} \\
\midrule
Noisy input        & 5.200 & 17.302 & 34.636 & 8.662 & 8.696 & 8.582 \\
\midrule
DFT \cite{smith2011spectral} & 2.710 & 8.376 & 16.425 & 5.948 & 13.361 & 27.704 \\
DLCT \cite{healy2016linear} & 2.682 & 8.306 & 16.233 & 5.815 & 13.099 & 27.040 \\
DFT \cite{hitzer2021introduction} & 2.333 & 7.032 & 13.466 & 3.405 & 3.456 & 3.393 \\
QLCT \cite{bahri2021discrete} & 2.320 & 7.005 & 13.446 & 3.369 & 3.414 & 3.356 \\
Savitzky--Golay \cite{savitzky1964} & 2.388 & 7.835 & 16.265 & 4.060 & 4.129 & 4.110 \\
RTS \cite{rauch1965} & 1.581 & 4.234 & 8.687 & 2.925 & 3.051 & 3.042 \\
\midrule
CA-DQ-DLCT & \bfseries 1.547 & \bfseries 3.425 & \bfseries 5.565 & \bfseries 2.782 & \bfseries 2.931 & \bfseries 2.894 \\
\bottomrule
\end{tabular}
\end{table}
\FloatBarrier

\subsubsection{Two-dimensional Denoising}\label{subsec:twoaxis}

We construct $256\times512$ pose fields from accelerating segments of
Mini Wheelbot encoder records \cite{hose2026wheelbot}, using two-axis
rotations and a $0.10$~m lever arm. Translation varies along the first
axis. Disjoint source-file pools provide $48$ fitting and $48$ test
fields.

Noise levels match Section~\ref{subsec:denoise1d}. The chirp grids are
$\{-0.2+0.0025r:r=0,\ldots,160\}$ and $\mathcal A_0$ for the two axes.
The same two-stage algorithm and fitting objectives used for sequences
apply over both axes.

The baselines follow the protocols in Section~\ref{subsec:denoise1d}
along both axes. DLCT and QLCT use rates $(0.0025,0.0025)$; DFT and
DFT use zero. Savitzky--Golay and RTS smoothing operate successively
along the two axes, without using the generating joint model. Fitting
selects an RTS noise-strength ratio of $10^{-2}$ for these fields.

CA-DQ-DLCT yields the lowest mean errors in Table~\ref{tab:twoaxis}.
Relative to DFT and QLCT, rotation RMSE decreases by
$59.7$--$85.2\%$ and position RMSE by $47.0$--$58.9\%$. At the
highest noise level, position RMSE is $3.952$~mm for RTS, compared
with $0.680$~mm for CA-DQ-DLCT. These results concern the evaluated
settings and motion models and demonstrate the performance of the
complete CA-DQ-DLCT pipeline. The individual effects of the
dual-quaternion representation, chirp adaptation, and position
refinement are not isolated by these comparisons.

\begin{table}[t]
\centering
\caption{Mean per-field rotation RMSE (degrees) and position RMSE
(mm) on $48$ test fields constructed from Mini Wheelbot trajectories.
Each field has size $256\times512$; errors are evaluated on its
central half along each axis. DLCT and QLCT use
$(\alpha_1,\alpha_2)=(0.0025,0.0025)$. Bold identifies column minima
before rounding.}
\label{tab:twoaxis}
\vspace{6pt}
\sisetup{table-format=2.3}
\small
\setlength{\tabcolsep}{3pt}
\begin{tabular}{l SSS SSS}
\toprule
& \multicolumn{3}{c}{Rotation RMSE (deg)} & \multicolumn{3}{c}{Position RMSE (mm)} \\
\cmidrule(lr){2-4}\cmidrule(lr){5-7}
Method & {$3^\circ$} & {$10^\circ$} & {$20^\circ$} & {$3^\circ$} & {$10^\circ$} & {$20^\circ$} \\
\midrule
Noisy input                               & 5.193 & 17.322 & 34.658 & 8.664 & 8.665 & 8.657 \\
\midrule
DFT \cite{smith2011spectral} & 2.093 & 7.877 & 16.094 & 3.950 & 5.460 & 8.890 \\
DLCT \cite{healy2016linear} & 2.066 & 7.791 & 15.927 & 3.887 & 5.248 & 8.402 \\
DFT \cite{hitzer2021introduction} & 1.573 & 6.213 & 12.759 & 1.656 & 1.657 & 1.653 \\
QLCT \cite{bahri2021discrete} & 1.579 & 6.221 & 12.767 & 1.569 & 1.570 & 1.565 \\
Savitzky--Golay \cite{savitzky1964} & 2.196 & 4.195 & 7.608 & 2.884 & 3.472 & 4.876 \\
RTS \cite{rauch1965} & 2.141 & 3.230 & 5.302 & 2.665 & 3.020 & 3.952 \\
\midrule
CA-DQ-DLCT & \bfseries 0.634 & \bfseries 1.375 & \bfseries 1.888 & \bfseries 0.832 & \bfseries 0.688 & \bfseries 0.680 \\
\bottomrule
\end{tabular}
\end{table}
\FloatBarrier

\section{Conclusion}\label{sec:conclusion}

This paper introduced the dual quaternion discrete linear canonical
transform (DQ-DLCT) in its one- and two-dimensional left-, right-, and
two-sided forms. A chirp--Fourier--chirp factorization yielded fast
forward and inverse algorithms, reducing the two-dimensional cost from
$\mathcal O(N_1^2N_2^2)$ to
$\mathcal O\!\left(N_1N_2(\log N_1+\log N_2)\right)$, and spatial and
frequency convolution theorems were established. A chirp-adaptive
DQ-DLCT method was applied to rigid-body pose denoising, and numerical
results confirmed the efficiency and effectiveness of the framework.

The approach combines the compact, singularity-free pose
representation of dual quaternions with the extra degrees of freedom
of canonical transforms, and the FFT-based algorithm makes this
combination practical. Its main limitations are the algebraic
complexity of noncommutative dual quaternion operations, the absence
of a closed-form rule for canonical-parameter selection, and
validation confined to measured sequences and constructed pose fields.

Both ingredient techniques are widely used---dual quaternions in
robot kinematics, motion interpolation, SLAM, and hand--eye
calibration, and linear canonical transforms in optics, radar and
sonar, time--frequency analysis, and colour image processing---while
their intersection remains an emerging direction. Future work will
address finite-precision stability, automatic parameter selection,
GPU implementations, benchmarking against mature pose filters, and
large-scale motion processing.

\section*{Acknowledgments}
The authors acknowledge partial funding from the University of Macau Development Foundation (UMDF) under Grants MYRG-GRG2024-00290-FST-UMDF and MYRG-GRG2025-00250-FST.

\section*{CRediT Authorship Contribution Statement}·         
\textbf{Yi Wu:} Conceptualization, Formal analysis, Writing -- original draft.
\textbf{Xuan Chen:} Software, Validation, Writing -- review \& editing.
\textbf{Kit Ian Kou:} Supervision, Writing -- review \& editing.

\section*{Declaration of Competing Interest}
The authors declare that they have no known competing financial interests or personal relationships that could have influenced the work reported in this paper.

\section*{Data Availability}
The Mini Wheelbot dataset used in this study is publicly available from the repository cited in the manuscript.

\bibliographystyle{unsrtnat}
\bibliography{ref}
\end{document}